\documentclass[15pt, a4paper]{article}
\usepackage{}
\usepackage{amsfonts}
\usepackage{mathbbold}
\usepackage{bbm}
\usepackage{mathrsfs}
\usepackage{latexsym}
\usepackage{graphicx}
\usepackage{amssymb}

\newtheorem{theorem}{Theorem}[section]
\newtheorem{lemma}[theorem]{Lemma}
\newtheorem{corollary}[theorem]{Corollary}

\title{{\Large \bf The uniform supertrees with the extremal spectral radius\thanks{Supported by National natural science foundation of China (NSFC)
(Nos. Nos. 12171222, 12101285, 12071411), Natural science foundation of Guangdong province (No. 2021A1515010254), Foundation of Lingnan Normal
University(ZL1923).}}}
\author{Guanglong Yu$^{a}$
~ Lin Sun$^{a}$\thanks{Corresponding authors, E-mail addresses:
yglong01@163.com (G. Yu), sunlin@lingnan.edu.cn (L. Sun).} ~  Hailiang Zhang$^b$ ~ Gang Li$^{c}$ ~
\\ ~ \\
{\footnotesize $^a$Department of Mathematics, Lingnan Normal
University,  Zhanjiang, Guangdong, 524048, P.R.China}\\
{\footnotesize $^b$Department of Mathematics, Taizhou University, Linhai, Zhejiang, 317000, P.R.China}\\
{\footnotesize $^c$College of Mathematics and Systems Science, Shandong University of Science and Technology,}\\
{\footnotesize  Qingdao, Shandong, 266590, P.R.China}}
\date{}

\begin{document}
%\openup 1.0\jot
\maketitle

\begin{abstract}
For a $hypergraph$ $\mathcal{G}=(V, E)$ consisting of a nonempty vertex set $V=V(\mathcal{G})$ and an edge set $E=E(\mathcal{G})$, its $adjacency$ $matrix$ $\mathcal {A}_{\mathcal{G}}=[(\mathcal {A}_{\mathcal{G}})_{ij}]$ is defined as
$(\mathcal {A}_{\mathcal{G}})_{ij}=\sum_{e\in E_{ij}}\frac{1}{|e| - 1}$, where $E_{ij} = \{e \in E : i, j \in e\}$.
The $spectral$ $radius$ of a hypergraph $\mathcal{G}$, denoted by $\rho(\mathcal {G})$, is the maximum modulus among all eigenvalues of $\mathcal {A}_{\mathcal{G}}$.
In this paper, among all $k$-uniform ($k\geq 3$) supertrees with fixed number of vertices, the supertrees with the maximum, the second maximum and the minimum spectral radius are completely determined, respectively.


\end{abstract}
\baselineskip 18.6pt

\section{Introduction}

\ \ \ \
In the past twenty years, some different connectivity hypermatrices (or tensors) had been defined and been developed to explore spectral hypergraph theory  \cite{Ban.CB}-\cite{KCKPZ.2}, \cite{DGAO}, \cite{KCHX}-\cite{HLBZ}, \cite{NgQZ}, \cite{LQIE}-\cite{S.SZ}, \cite{YQY, YQY2}. Using different hypermatrices for general hypergraphs, many interesting spectral properties have been studied that many properties of spectral graph theory have been
extended to spectral hypergraph theory. A lot of interesting results have emerged and the spectra of hypergraphs have been further studied \cite{CJDA, YFTP, SFSH, HGBZ,  LLKB, LLSM, COQY, PXLW}, \cite{GYSWZ}-\cite{JZLS}. In \cite{A.Ban}, A. Banerjee introduced an adjacency matrix and use its spectrum so that some spectral and structural
properties of hypergraphs are revealed. In this paper, we go on studying the spectra of hypergraphs according to the adjacency matrix introduced in \cite{A.Ban}.

Now we recall some notations and definitions related to hypergraphs.
For a set $S$, we denote by $|S|$ its cardinality. A $hypergraph$ $\mathcal{G}=(V, E)$ consists of a nonempty vertex set $V=V(\mathcal{G})$ and an edge set $E=E(\mathcal{G})$, where each edge $e\in E(\mathcal{G})$ is a subset of $V(\mathcal{G})$ containing at least two vertices. The cardinality $n=|V(\mathcal{G})|$ is called the order; $m=|E(\mathcal{G})|$ is called the edge number of hypergraph $\mathcal{G}$. Denote by $t$-set a set with size (cardinality) $t$. We say that a hypergraph $\mathcal{G}$ is $uniform$ if its every edge has the same size, and call it $k$-$uniform$ if its every edge has size $k$ (i.e. every edge is a $k$-subset). It is known that a $2$-uniform graph is always called a ordinary graph or graph for short.

For a hypergraph $\mathcal{G}$, we define $\mathcal{G}-e$ ($\mathcal{G}+e$)
to be the hypergraph obtained from $\mathcal{G}$ by deleting the edge $e\in
 E(\mathcal{G})$ (by adding a new edge $e$ if $e\notin
 E(\mathcal{G})$); for an edge subset $B\subseteq E(\mathcal{G})$, we define $\mathcal{G}-B$
to be the hypergraph obtained from $\mathcal{G}$ by deleting each edge $e\in
 B$; for a vertex subset $S\subseteq V(\mathcal{G})$, we define $\mathcal{G}-S$ to be the hypergraph obtained from $\mathcal{G}$ by deleting all the vertices in $S$ and deleting the edges incident with any vertex in $S$. For two $k$-uniform hypergraphs $\mathcal{G}_{1}=(V_{1}, E_{1})$ and $\mathcal{G}_{2}=(V_{2}, E_{2})$, we say the two graphs are $isomorphic$ if there is a bijection $f$ from $V_{1}$ to $V_{2}$, and there is a bijection $g$ from $E_{1}$ to $E_{2}$ that maps each edge $\{v_{1}$, $v_{2}$, $\ldots$, $v_{k}\}$ to $\{f(v_{1})$, $f(v_{2})$, $\ldots$, $f(v_{k})\}$.

In a hypergraph, two vertices are said to be $adjacent$ if
both of them are contained in an edge. Two edges are said to be $adjacent$
if their intersection is not empty. An edge $e$ is said to be $incident$ with a vertex $v$ if
$v\in e$. The $neighbor$ $set$ of vertex $v$ in hypergraph $\mathcal{G}$, denoted by $N_{\mathcal{G}}(v)$, is the set of vertices adjacent to $v$ in $\mathcal{G}$. The $degree$ of a vertex $v$ in $\mathcal{G}$, denoted by $deg_{\mathcal{G}}(v)$ (or $deg(v)$ for short), is the number of the edges incident with $v$. For a hypergraph $\mathcal{G}$, among all of its vertices, we denote by $\Delta(\mathcal{G})$ (or $\Delta$ for short) the $maximal$ $degree$, and denote by $\delta(\mathcal{G})$ (or $\delta$ for short) the $minimal$ $degree$ respectively. A vertex of degree $1$ is called a $pendant$ $vertex$. A $pendant$ $edge$ is an edge with at most one vertex of degree more than one and other vertices in this edge being all pendant vertices.

In a hypergraph, a $hyperpath$ of length $q$ ($q$-$hyperpath$) is defined to be an alternating sequence
of vertices and edges $v_{1}e_{1}v_{2}e_{2}\cdots v_{q}e_{q}v_{q+1}$ such that
(1) $v_{1}$, $v_{2}$, $\ldots$, $v_{q+1}$ are all distinct vertices;
(2) $e_{1}$, $e_{2}$, $\ldots$, $e_{q}$ are all distinct edges;
(3) $v_{i}$, $v_{i+1}\in e_{i}$ for $i = 1$, $2$, $\ldots$, $q$;
(4) $e_{i}\cap e_{i+1}=v_{i+1}$ for $i = 1$, $2$, $\ldots$, $q-1$; (5) $e_{i}\cap e_{j}=\emptyset$ if $|i-j|\geq 2$.
If there is no discrimination, a hyperpath is sometimes written as $e_{1}e_{2}\cdots e_{q-1}e_{q}$, $e_{1}v_{2}e_{2}\cdots v_{q}e_{q}$ or $v_{1}e_{1}v_{2}e_{2}\cdots v_{q}e_{q}$. A $hypercycle$ of length $q$ ($q$-$hypercycle$) $v_{1}e_{1}v_{2}e_{2}\cdots v_{q-1}e_{q-1}v_{q}e_{q}v_{1}$ is obtained from a hyperpath $v_{1}e_{1}v_{2}e_{2}\cdots v_{q-1}e_{q-1}v_{q}$ by adding a new edge $e_{q}$ between $v_{1}$ and $v_{q}$ where $e_{q}\cap e_{1}= \{v_{1}\}$, $e_{q}\cap e_{q-1}= \{v_{q}\}$, $e_{q}\cap e_{j}=\emptyset$ if $j\neq1, q-1$ and $|q-j|\geq 2$. The length of a hyperpath $P$ (or a hypercycle $C$), denoted by $L(P)$ (or $L(C)$), is the number of the edges in $P$ (or $C$). A hypergraph $\mathcal{G}$
is connected if there exists a hyperpath from $v$ to $u$ for all $v, u \in V$,
and $\mathcal{G}$ is called $acyclic$ if it contains no hypercycle.

Recall that a tree is an ordinary graph which is 2-uniform, connected and acyclic. A $supertree$ is similarly defined to be a hypergraph which is both connected and acyclic.
Clearly, in a supertree, its each pair of the
edges have at most one common vertex. Therefore, the edge number of a $k$-uniform supertree of order $n$ is $m=\frac{n-1}{k-1}$.

Let $G = (V, E)$ be an ordinary graph (2-uniform). For every $k \geq 3$,
the kth power of $G$, denoted by $G^{k}= (V^{k} , E^{k} )$, is defined as the $k$-uniform hypergraph with
the edge set $E^{k}= \{e \cup \{v_{e_{1}}$, $v_{e_{2}}$, $\ldots$, $v_{e_{k-2}}\}:e\in E\}$ and the vertex set $V^{k}= V \cup
(\cup_{e\in E} \{v_{e_{1}}$, $v_{e_{2}}$, $\ldots$, $v_{e_{k-2}}\})$, where $V \cap
(\cup_{e\in E} \{v_{e_{1}}$, $v_{e_{2}}$, $\ldots$, $v_{e_{k-2}}\})=\emptyset$, $\{v_{e_{1}}$, $v_{e_{2}}$, $\ldots$, $v_{e_{k-2}}\}\cap \{v_{f_{1}}$, $v_{f_{2}}$, $\ldots$, $v_{f_{k-2}}\}=\emptyset$ for $e\neq f$, $e,f\in E$. The kth power of an ordinary tree is called a $hypertree$. Obviously, a hypertree is a supertree.

Denote by $\mathcal{P}(n, k)$ the $k$-uniform hyperpath of order $n$. A $k$-uniform $superstar$ of order $n$, denoted by $\mathcal{S}^{\ast}(n, k)$ (see Fig. 1.1), is a supertree in which all edges intersect at just one common vertex. A $k$-uniform $double$ $hyperstar$ of order $n$, denote by $\mathcal{S}(n, k; l_{1}, l_{2})$ where $l_{1}, l_{2}\geq 1$ (see Fig. 1.1), is a supertree obtained by attaching $l_{1}$ pendant edges at vertex $u_{1}$ of an edge $e$, and attaching $l_{2}$ pendant edges at the other vertex $u_{2}$ of edge $e$, where $u_{1}\neq u_{2}$.

\setlength{\unitlength}{0.6pt}
\begin{center}
\begin{picture}(669,165)
\put(111,46){\circle*{4}}
\qbezier(3,93)(0,61)(111,46)
\qbezier(3,93)(31,124)(111,46)
\qbezier(111,46)(42,143)(69,155)
\qbezier(111,46)(101,165)(69,155)
\qbezier(111,46)(216,73)(211,109)
\qbezier(111,46)(188,128)(211,109)
\put(123,138){\circle*{4}}
\put(136,133){\circle*{4}}
\put(150,126){\circle*{4}}
\put(450,62){\circle*{4}}
\qbezier(391,140)(361,122)(450,62)
\qbezier(391,140)(409,157)(450,62)
\qbezier(450,62)(369,43)(352,61)
\qbezier(352,61)(350,92)(450,62)
\put(570,62){\circle*{4}}
\qbezier(450,62)(512,43)(570,62)
\qbezier(450,62)(515,88)(570,62)
\qbezier(570,62)(658,44)(669,64)
\qbezier(570,62)(665,85)(669,64)
\qbezier(570,62)(587,132)(629,137)
\qbezier(629,137)(646,115)(570,62)
\put(104,24){$\mathcal{S}^{\ast}(n, k)$}
\put(481,24){$\mathcal{S}(n, k; 2, 2)$}
\put(442,49){$u_{1}$}\put(567,49){$u_{2}$}\put(507,79){$e$}
\put(225,-9){Fig. 1.1. $\mathcal{S}^{\ast}(n, k)$ and $\mathcal{S}(n; 2, 2)$}
\end{picture}
\end{center}

Let $E_{ij} = \{e \in E : i, j \in e\}$. The $adjacency$ $matrix$ $\mathcal {A}_{\mathcal{G}}=[(\mathcal {A}_{\mathcal{G}})_{ij}]$ of a hypergraph $\mathcal{G}$ is defined as
$$(\mathcal {A}_{\mathcal{G}})_{ij}=\sum_{e\in E_{ij}}\frac{1}{|e| - 1}.$$ It is easy to find that $\mathcal {A}_{\mathcal{G}}$ is symmetric if there is no requirement for direction on hypergraph $\mathcal{G}$, and find that $\mathcal {A}_{\mathcal{G}}$ is very convenient to be used to investigate the spectum of a hypergraph even without the requirement for edge uniformity. The $spectral$ $radius$ $\rho(\mathcal {G})$ of a hypergraph $\mathcal{G}$ is defined to be the spectral radius $\rho(\mathcal {A}_{\mathcal{G}})$, which is the maximum modulus among all eigenvalues of $\mathcal {A}_{\mathcal{G}}$.
In spectral theory of hypergraphs, the spectral radius is an index that attracts much attention due to its fine properties \cite{KCKPZ, YFTP, SFSH, LSQ, HLBZ, LLSM, COQY, LSLK, YQY, CLYXZ, JZLG}.

We assume that the hypergraphs throughout
this paper are simple, i.e. $e_{i} \neq e_{j}$ if $i \neq j$, and assume the hypergraphs throughout
this paper are undirected. In this paper, among all $k$-uniform ($k\geq 3$) supertrees with fixed number of vertices, the supertrees with the maximum, the second maimum and the minimum spectral radius are completely determined respectively, getting the following result:

\begin{theorem}\label{th01.01} %------
Let $\mathcal{G}$ be a $k$-uniform  ($k\geq 3$) supertree of order $n$. Then $\rho(\mathcal{G})\leq\rho(\mathcal{S}^{\ast}(n, k))$ with equality if and only if $\mathcal{G}\cong \mathcal{S}^{\ast}(n, k)$.

\end{theorem}

\begin{theorem}\label{th01.02} %------
Let $\mathcal{G}$ be a $k$-uniform  ($k\geq 3$) supertree of order $n$ and with $m(\mathcal{G})\geq 3$ satisfying that $\mathcal{G}\ncong \mathcal{S}^{\ast}(n, k)$. Then $\rho(\mathcal{G})\leq\rho(\mathcal{S}(n, k; \frac{n-1}{k-1}-2, 1))$ with equality if and only if $\mathcal{G}\cong \mathcal{S}(n, k; \frac{n-1}{k-1}-2, 1)$.

\end{theorem}

\begin{theorem}\label{th01.03} %------
Let $\mathcal{G}$ be a $k$-uniform ($k\geq 3$) supertree of order $n$. Then $\rho(\mathcal{P}(n, k))\leq \rho(\mathcal{G})$ with equality if and only if $\mathcal{G}\cong \mathcal{P}(n, k)$.

\end{theorem}

\begin{corollary} \label{cl01,04}
Suppose $T^{k}$ ($k\geq 3$) of order $n$ is the kth power of ordinary tree $T$. Then

(1) $\rho(T^{k})\leq \mathcal{S}^{\ast}(n, k)$ with equality if and only if $T^{k}\cong \mathcal{S}^{\ast}(n, k)$.

(2)
$\rho(\mathcal{P}(n, k))\leq \rho(T^{k})$ with equality if and only if $T^{k}\cong \mathcal{P}(n, k)$.
\end{corollary}

The layout of this paper is as follows: section 2 introduces some basic knowledge and working lemmas; section 3 represents our results.

\section{Preliminary}

\ \ \ \ \ For the requirements in the narrations afterward, we need some prepares. For a hypergraph $\mathcal{G}$ with vertex set $\{v_{1}$, $v_{2}$, $\ldots$, $v_{n}\}$, a vector on $\mathcal{G}$ is a vector $X=(x_{v_1}, x_{v_2}, \ldots, x_{v_n})^T \in R^n$ that
 entry $x_{v_i}$ is mapped to vertex $v_i$ for $i\leq i\leq n$.

 From \cite{OP}, by the famous Perron-Frobenius theorem, for $\mathcal {A}_{\mathcal{G}}$ of a connected uniform hypergraph $\mathcal{G}$ of order $n$, we know that there is unique one positive eigenvector $X=(x_{v_{1}}$, $x_{v_{2}}$, $\ldots$, $x_{v_{n}})^T \in R^{n}_{++}$ ($R^{n}_{++}$ means the set of positive real vectors of dimension $n$) corresponding to $\rho(\mathcal{G})$, where $\sum^{n}_{i=1}x^{2}_{v_{i}}= 1$ and each entry $x_{v_i}$ is mapped to each vertex $v_i$ for $i\leq i\leq n$. We call such an eigenvector $X$
the $principal$ $eigenvector$ of $\mathcal{G}$.

Let $A$ be an irreducible nonnegative $n \times n$ real matrix (with every entry being real number) with spectral radius $\rho(A)$. The following extremal representation (Rayleigh quotient) will be useful:
$$\rho(A)=\max_{X\in R^{n}, X\neq0}\frac{X^{T}AX}{X^{T}X},$$ and if a vector $X$ satisfies that $\frac{X^{T}AX}{X^{T}X}=\rho(A)$, then $AX=\rho(A)X$.

\begin{lemma} \label{le3,05}
Let $A$ be an irreducible nonnegative square real matrix with order $n$
and spectral radius $\rho$, $Y \in (R^{n}_{+}\setminus \{0\}^{n})$ be a nonnegative vector ($R^{n}_{+}$ means the set of nonnegative real vectors of dimension $n$, $\{0\}^{n}=\{(0, 0, \ldots, 0)^{T}\}$). If $AY\geq \rho Y$, then
$AY= \rho Y$.
\end{lemma}

\begin{proof}
Note the relation between the spectral radius  and Rayleigh quotient for an irreducible nonnegative square real matrix. It follows that $\frac{Y^{T}AY}{Y^{T}Y}= \rho$, and $AY=\rho Y$. Thus the result follows. This completes the proof. \ \ \ \ \ $\Box$
\end{proof}

\begin{lemma}{\bf \cite{GYJSH}} \label{le3,04}
Let $A$ be an irreducible nonnegative square symmetric real matrix with order $n$
and spectral radius $\rho$, $Y \in (R^{n}_{+}\setminus \{0\}^{n})$ be a nonnegative vector. If there exists $r\in R_{+}$
such that $AY\leq rY$, then
$\rho\leq r$. Similarly, if there exists $r\in R_{+}$
such that $AY\geq rY$, then
$\rho\geq r$.
\end{lemma}

\section{Main results}

~~~~Let $X$ be an eigenvector of a connected $k$-uniform hypergraph $G$. For
the simplicity, we let $\displaystyle x_{e} =\sum_{i<j, v_{i}, v_{j}\in e} x_{v_{i}}x_{v_{j}}$ for an edge $e=\{v_{1}$, $v_{2}$, $\ldots$, $v_{k}\}$.

\begin{lemma}\label{le03.02} %------
Let $e_{1}=\{v_{1,1}$, $v_{1,2}$, $\ldots$, $v_{1,k}\}$, $e_{2}=\{v_{2,1}$, $v_{2,2}$, $\ldots$, $v_{2,k}\}$, $\ldots$, $e_{j}=\{v_{j,1}$, $v_{j,2}$, $\ldots$, $v_{j,k}\}$ be some edges in a connected $k$-uniform hypergraph $\mathcal{G}$; $v_{u,1}$, $v_{u,2}$, $\ldots$, $v_{u,t}$ be vertices in $\mathcal{G}$ that $t< k$. For $1\leq i\leq j$, $\{v_{u,1}, v_{u,2}, \ldots, v_{u,t}\}\nsubseteq e_{i}$, $e^{'}_{i}=(e_{i}\setminus \{v_{i,1}, v_{i,2}, \ldots, v_{i,t}\})\cup\{v_{u,1}, v_{u,2}, \ldots, v_{u,t}\}$ satisfying that $e^{'}_{i}\notin E(\mathcal{G})$. Let $\mathcal{G}^{'}=\mathcal{G}-\sum e_{i}+\sum e^{'}_{i}$. If in the principal eigenvector $X$ of $\mathcal{G}$, for $1\leq i\leq j$, $x_{v_{i,1}}\leq x_{v_{u,1}}$, $x_{v_{i,2}}\leq x_{v_{u,2}}$, $\ldots$, and $x_{v_{i,t}}\leq x_{v_{u,t}}$, then $\rho(\mathcal{G}^{'})>\rho(\mathcal{G})$.
\end{lemma}

\begin{proof}
Note that $X^{T}(\mathcal {A}_{\mathcal{G^{'}}}-\mathcal {A}_{\mathcal{G}})X=\frac{2}{k-1}\sum(x_{e^{'}_{i}}-x_{e_{i}})\geq 0$. It follows that $\rho(\mathcal{G}^{'})\geq\rho(\mathcal{G})$.
If $\rho(\mathcal{G}^{'})=\rho(\mathcal{G})$, then $\rho(\mathcal{G}^{'})=X^{T}\mathcal {A}_{\mathcal{\mathcal{G}^{'}}}X=X^{T}\mathcal {A}_{\mathcal{G}}X=\rho(\mathcal{G})$. It follows that $\mathcal {A}_{\mathcal{\mathcal{G}^{'}}}X=\rho(\mathcal{G}^{'})X=\rho(\mathcal{G})X=\mathcal {A}_{\mathcal{G}}X$. Without loss of generality, suppose $v_{u,1}\notin e_{1}$. Then
$(\mathcal {A}_{\mathcal{\mathcal{G}^{'}}}X)_{v_{u,1}}-(\mathcal {A}_{\mathcal{G}}X)_{v_{u,1}}\geq \frac{x_{v_{1,k}}}{k-1}>0$, which contradicts $\mathcal {A}_{\mathcal{G^{'}}}X=\mathcal {A}_{\mathcal{G}}X$. Consequently, it follows that $\rho(\mathcal{G}^{'})>\rho(\mathcal{G})$.
This completes the proof. \ \ \ \ \ $\Box$
\end{proof}

\begin{Proof}
Suppose $\mathcal{T}$ is a $k$-uniform supertree of order $n$ satisfying that $\rho(\mathcal{T})=\max\{\rho(\mathcal{G}):\, \mathcal{G}$ is a $k$-uniform supertree of order $n\}$.
Let $X$ be
the principal eigenvector of $\mathcal{T}$ and $x_{u}=\max\{x_{v}: v\in V(\mathcal{T})\}$. Suppose that $\mathcal{T}\ncong \mathcal{S}^{\ast}(n, k)$. Then in $\mathcal{T}$, there exist edges not incident with vertex $u$. Suppose $e^{'}$ is not incident with vertex $u$. Note that $\mathcal{T}$ is connected. Then there is a hyperpath $P=ue_{1}v_{1}e_{2}v_{2}\cdots v_{t}e^{'}$ from $u$ to $e^{'}$. Let $e^{'}_{2}=(e_{2}\setminus \{v_{1}\})\cup \{u\}$, and $\mathcal{T}^{'}=\mathcal{T}-e_{2}+e^{'}_{2}$. Then by Lemma \ref{le03.02}, it follows that $\rho(\mathcal{T}^{'})>\rho(\mathcal{T})$, which contradicts the maximality of $\rho(\mathcal{T})$. Hence, it follows that $\mathcal{T}\cong \mathcal{S}^{\ast}(n, k)$. This completes the proof. \ \ \ \ \ $\Box$
\end{Proof}

\begin{Pro}
Let $\Lambda=\{\mathcal{G}:\, \mathcal{G}$ be a $k$-uniform supertree of order $n$ and $\mathcal{G}\ncong \mathcal{S}^{\ast}(n, k)\}$. Suppose $\mathcal{T}\in \Lambda$ satisfies that $\rho(\mathcal{T})=\max\{\rho(\mathcal{G}):\, \mathcal{G}\in\Lambda\}$.
Let $X$ be
the principal eigenvector of $\mathcal{T}$ and $x_{u}=\max\{x_{v}: v\in V(\mathcal{T})\}$. Note that $\mathcal{T}\ncong \mathcal{S}^{\ast}(n, k)$, and $\mathcal{T}$ is connected. Then in $\mathcal{T}$, there exist a hyperpath $P=ue_{1}v_{1}e_{2}$.

{\bf Assertion} Except edges $e_{1}$, $e_{2}$, any one of other edges is incident with vertex $u$. Otherwise, suppose one edge $e_{t}$ is not incident with vertex $u$. Note that $\mathcal{T}$ is connected. Then there is a hyperpath $\mathcal{P}=ue_{a_{1}}v_{a_{1}}e_{a_{2}}v_{a_{2}}\cdots v_{a_{t}}e_{t}$ from $u$ to $e_{t}$.

{\bf Claim 1} If $e_{1}\in E(\mathcal{P})$, then $e_{a_{1}}=e_{1}$. Otherwise, suppose $e_{a_{i}}=e_{1}$ where $1< i\leq t$. Then $ue_{a_{1}}v_{a_{1}}e_{a_{2}}v_{a_{2}}\cdots v_{a_{i}}e_{a_{i}}u$ contains cycle, which contradicts that $\mathcal{T}$ is a supertree.

In the same way, we get the following Claim 2.

{\bf Claim 2} If $e_{2}\in E(\mathcal{P})$, then $e_{a_{1}}=e_{1}$, $e_{a_{2}}=e_{2}$.

Therefore, there 3 cases to consider, which are: (1) $e_{a_{1}}=e_{1}$, $e_{2}\notin E(\mathcal{P})$; (2) $e_{a_{1}}=e_{1}$, $e_{a_{2}}=e_{2}$; (3) $e_{1}\notin E(\mathcal{P})$, $e_{2}\notin E(\mathcal{P})$. For the case that $e_{a_{1}}=e_{1}$, $e_{2}\notin E(\mathcal{P})$, let $e^{'}_{a_{2}}=(e_{a_{2}}\setminus \{v_{a_{1}}\})\cup \{u\}$, and $\mathcal{T}^{'}=\mathcal{T}-e_{a_{2}}+e^{'}_{a_{2}}$, where $\mathcal{T}^{'}\in \Lambda$. Then by Lemma \ref{le03.02}, it follows that $\rho(\mathcal{T}^{'})>\rho(\mathcal{T})$, which contradicts the maximality of $\rho(\mathcal{T})$. In the same way, for the cases (2) and (3), we can get a supertree $\mathcal{T}^{'}$ where $\mathcal{T}^{'}\in \Lambda$, such that $\rho(\mathcal{T}^{'})>\rho(\mathcal{T})$ which contradicts the maximality of $\rho(\mathcal{T})$. Thus, our assertion holds.

From the above assertion, it follows that $\mathcal{T}\cong \mathcal{S}(n, k; \frac{n-1}{k-1}-2, 1)$. This completes the proof. \ \ \ \ \ $\Box$
\end{Pro}

\begin{lemma}{\bf \cite{H.M}} \label{le3,06}
Let $A$ be an irreducible nonnegative square matrix with order $n$
and spectral radius $\rho$. Let $s^{A}_{i}$ be the ith row sum, $s_{A}=\min\{s^{A}_{i}: 1\leq i\leq n\}$, and $S_{A}=\max\{s^{A}_{i}: 1\leq i\leq n\}$. Then $s_{A}\leq\rho\leq S_{A}$ with either one equality if and only if $A$ is regular (all of the row sums of $A$ are equal).
\end{lemma}

From Lemma \ref{le3,06}, combining with hypergraph, we can get the following corollary naturally.

\begin{corollary} \label{le3,07}
For a connected hypergraph $\mathcal{G}$, we have $\delta\leq\rho(\mathcal{G})\leq \Delta$ with either one equality if and only if $\mathcal{G}$ is regular, where $\delta$ is the minimum degree, $\Delta$ is the maximum degree.
\end{corollary}

Using Lemma \ref{le3,06}, we can get an improvement for Lemma \ref{le3,04}.

\begin{lemma} \label{le3,04,01}
Let $A$ be an irreducible nonnegative square symmetric real matrix with order $n$
and spectral radius $\rho$, $y\in R^{n}_{++}$ be a positive vector. If there exists $r\in R_{+}$ such that $Ay\leq ry$, then
$\rho\leq r$ with equality if and only if $Ay= ry$. Similarly, if there exists $r\in R_{+}$ such that $Ay\geq ry$, then
$\rho\geq r$ with equality if and only if $Ay= ry$.
\end{lemma}

\begin{proof}
Using Lemma \ref{le3,04} gets that $\rho\leq r$ if $Ay\leq ry$; $\rho\geq r$ if $Ay\geq ry$. Next we prove the conclusion for $\rho= r$.

We first prove the conclusion that $\rho= r$ if and only if $Ay= ry$ under the condition $Ay\leq ry$. Suppose $y=(y_{1}, y_{2}, \ldots, y_{n})^{T}$. Let $B = \left(\begin{array}{cccccc}
\frac{1}{y_{1}} & 0 & 0 & \cdots & 0\\
0 & \frac{1}{y_{2}} & 0 & \vdots &  0 \\
\vdots \\
 0 & 0 & \vdots & 0 & \frac{1}{y_{n}}\\
\end{array}\right)A\left(\begin{array}{cccccc}
y_{1} & 0 & 0 & \cdots & 0\\
0 & y_{2} & 0 & \vdots &  0 \\
\vdots \\
 0 & 0 & \vdots & 0 & y_{n}\\
\end{array}\right)$. Denote by $\rho(B)$ the spectral radius of $B$. Note that the eigenvalues of $B$ are the same to the eigenvalues of $A$; $\rho= r$ means $\rho(B)= r$. Note that $Ay\leq ry$ means $S_{B}\leq r$; $\rho(B)= r$ means all of the row sums of $B$ equals $r$ by Lemma \ref{le3,06}, which implies that $Ay= ry$. As a result, it follows that under the condition that $Ay\leq ry$, if $\rho= r$, then $Ay= ry$. Conversely, if $Ay= ry$, then all of the row sums of $B$ equals $r$, and then $\rho(B)= r= \rho$.

In the same way, we get that $\rho= r$ if and only if $Ay= ry$ under the condition that $Ay\geq ry$. This completes the proof. \ \ \ \ \ $\Box$
\end{proof}

\begin{lemma} \label{le3,07,01}
(1) Suppose $c>0$, $d>0$, $a-c> 0$, $b-d>0$. If $\frac{a}{b}\geq\frac{c}{d}$, then $\frac{a-c}{b-d}\geq\frac{a}{b}$ with equality if and only if $\frac{a}{b}=\frac{c}{d}$. Moreover, if $\frac{a}{b}>\frac{c}{d}$, then $\frac{a-c}{b-d}>\frac{a}{b}$.

(2) Suppose $c>0$, $d>0$, $a-c> 0$, $b-d>0$. If $\frac{a}{b}\geq\frac{c}{d}$, then $\frac{a+c}{b+d}\leq\frac{a}{b}$ with equality if and only if $\frac{a}{b}=\frac{c}{d}$. Moreover, if $\frac{a}{b}>\frac{c}{d}$, then $\frac{a+c}{b+d}<\frac{a}{b}$.

(3) Suppose $\frac{a}{b}\geq1$, $b>c> 0$. Then $\frac{a-c}{b-c}\geq\frac{a}{b}$.
\end{lemma}

\begin{proof}
(1) From $\frac{a}{b}\geq\frac{c}{d}$, it follows that $ab-bc\geq ab-ad$, which induces $\frac{a-c}{b-d}\geq\frac{a}{b}$. In the same way, we get that $\frac{a-c}{b-d}>\frac{a}{b}$ if $\frac{a}{b}>\frac{c}{d}$. Then (1) follows.

(2) is proved as (1). (3) is a corollary following from (1).
This completes the proof. \ \ \ \ \ $\Box$
\end{proof}

\

\setlength{\unitlength}{0.6pt}
\begin{center}
\begin{picture}(578,116)
\qbezier(0,60)(0,72)(28,81)\qbezier(28,81)(57,90)(98,90)\qbezier(98,90)(138,90)(167,81)\qbezier(167,81)(196,72)(196,60)\qbezier(196,60)(196,47)(167,38)
\qbezier(167,38)(138,30)(98,30)\qbezier(98,30)(57,30)(28,38)\qbezier(28,38)(0,47)(0,60)
\qbezier(132,60)(132,74)(176,85)\qbezier(176,85)(221,95)(285,95)\qbezier(285,95)(348,95)(393,85)\qbezier(393,85)(438,74)(438,60)\qbezier(438,60)(438,45)(393,34)
\qbezier(393,34)(348,24)(285,24)\qbezier(285,25)(221,25)(176,34)\qbezier(176,34)(131,45)(132,60)
\qbezier(373,60)(373,72)(402,80)\qbezier(402,80)(432,89)(475,89)\qbezier(475,89)(517,89)(547,80)\qbezier(547,80)(577,72)(577,60)\qbezier(577,60)(577,47)(547,39)
\qbezier(547,39)(517,30)(475,30)\qbezier(475,31)(432,31)(402,39)\qbezier(402,39)(372,47)(373,60)
\put(310,65){\circle*{4}}
\put(299,65){\circle*{4}}
\put(164,65){\circle*{4}}
\put(157,51){$v_{1}$}
\put(289,65){\circle*{4}}
\put(252,65){\circle*{4}}
\put(221,65){\circle*{4}}
\put(406,66){\circle*{4}}
\put(399,51){$v_{k}$}
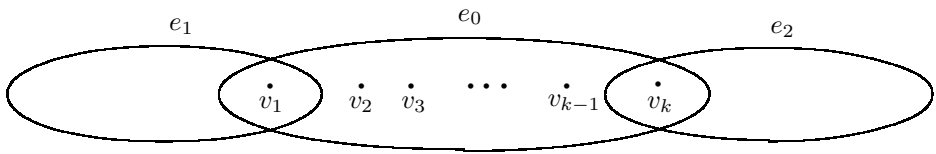
\put(195,-9){Fig. 3.1. $e_{0}$, $e_{1}$, $e_{2}$ in $\mathcal{G}$}
\put(101,100){$e_{1}$}
\put(281,106){$e_{0}$}
\put(476,98){$e_{2}$}
\put(213,51){$v_{2}$}
\put(246,51){$v_{3}$}
\put(349,65){\circle*{4}}
\put(337,52){$v_{k-1}$}
\end{picture}
\end{center}

\begin{lemma} \label{le3,08}
Let $\mathcal{G}$ be a hypergraph with spectral radius $\rho$, $e_{0}$, $e_{1}$, $e_{2}$ be three edges in $\mathcal{G}$ with $e_{0}=\{v_{1}$, $v_{2}$, $\ldots$, $v_{k-1}$, $v_{k}\}$, satisfying that $deg_{\mathcal{G}}(v_{2})=deg_{\mathcal{G}}(v_{3})=\cdots deg_{\mathcal{G}}(v_{k-1})=1$ ($k\geq 3$), $e_{1}\cap e_{0}=\{v_{1}\}$, $e_{2}\cap e_{0}=\{v_{k}\}$ (see Fig. 3.1.). Let $X$ be
the $principal$ $eigenvector$ of hypergraph $\mathcal{G}$. Then $x_{v_{2}}=x_{v_{3}}=\cdots=x_{v_{k-1}}=\frac{x_{v_{1}}+x_{v_{k}}}{(k-1)\rho-(k-3)}<\min\{x_{v_{1}}, x_{v_{k}}\}$.
\end{lemma}

\begin{proof}
For $2\leq i\leq k-1$, we prove $x_{v_{i}}<\min\{x_{v_{1}}, x_{v_{k}}\}$ by contradiction. Suppose that $\min\{x_{v_{1}}, x_{v_{k}}\}=x_{v_{1}}$, and $x_{v_{z}}\geq\min\{x_{v_{1}}, x_{v_{k}}\}$ for some $2\leq z\leq k-1$. Let $e^{'}_{1}=(e_{1}\setminus \{v_{1}\})\cup\{v_{z}\}$ and $\mathcal{G}_{1}=\mathcal{G}-e_{1}+e^{'}_{1}$. Using Lemma \ref{le03.02} gets $\rho(\mathcal{G}_{1})>\rho(\mathcal{G})$. But it contradicts $\rho(\mathcal{G}_{1})=\rho(\mathcal{G})$ because $\mathcal{G}_{1}\cong \mathcal{G}$. As a result, for $2\leq i\leq k-1$, it follows that $x_{v_{i}}<\min\{x_{v_{1}}, x_{v_{k}}\}$.

Note that $\rho x_{v_{2}}=\frac{1}{k-1}(x_{v_{1}}+x_{v_{3}}+\sum^{k}_{i=4}x_{v_{i}})$, $\rho x_{v_{3}}=\frac{1}{k-1}(x_{v_{1}}+x_{v_{2}}+\sum^{k}_{i=4}x_{v_{i}})$. It follows that $(\rho+\frac{1}{k-1})(x_{v_{2}}-x_{v_{3}})=0$. Note that $\rho> 1$ by Corollary \ref{le3,07}. Then we get $x_{v_{2}}=x_{v_{3}}$. Proceeding like this, we get that $x_{v_{2}}=x_{v_{3}}=\cdots=x_{v_{k-1}}$. Thus from $\rho x_{v_{2}}=\frac{1}{k-1}((k-3)x_{v_{3}}+x_{v_{1}}+x_{v_{k}})$, it follows that $x_{v_{2}}= \frac{x_{v_{1}}+x_{v_{k}}}{(k-1)\rho-(k-3)}$. Thus the result follows.
This completes the proof. \ \ \ \ \ $\Box$
\end{proof}

Simillar to Lemma \ref{le3,08}, we get the following Lemma \ref{le3,08,02}.

\begin{lemma} \label{le3,08,02}
Let $\mathcal{G}$ be a hypergraph with spectral radius $\rho$, $e=\{u$, $v_{1}$, $v_{2}$, $\ldots$, $v_{k-1}\}$ be a pendant edge in $\mathcal{G}$ ($k\geq 2$), where $deg_{\mathcal{G}}(u)\geq 2$. Then in the principal eigenvector $X$ of $\mathcal{G}$, $x_{v_{1}}=x_{v_{2}}=\cdots=x_{v_{k-1}}=\frac{x_{u}}{(k-1)\rho-(k-2)}<x_{u}$.
\end{lemma}

\setlength{\unitlength}{0.65pt}
\begin{center}
\begin{picture}(641,431)
\put(250,-9){Fig. 3.2. $\mathcal{G}_{0}$, $\mathcal{G}_{1}$, $\mathcal{G}_{2}$}
\put(331,404){\circle*{4}}
\put(387,405){\circle*{4}}
\put(267,404){\circle*{4}}
\put(290,273){\circle*{4}}
\put(352,273){\circle*{4}}
\put(411,57){\circle*{4}}
\put(465,57){\circle*{4}}
\put(495,178){\circle*{4}}
\put(548,177){\circle*{4}}
\put(437,178){\circle*{4}}
\put(327,391){$v_{1}$}
\put(293,402){$e_{1}$}
\put(352,403){$e_{2}$}
\put(265,390){$v_{2}$}
\put(378,391){$v_{3}$}
\put(404,402){$e_{3}$}
\put(235,402){$e_{4}$}
\put(430,402){\circle*{4}}
\put(441,402){\circle*{4}}
\put(452,402){\circle*{4}}
\put(196,404){\circle*{4}}
\put(207,404){\circle*{4}}
\put(218,404){\circle*{4}}
\put(402,274){\circle*{4}}
\put(413,274){\circle*{4}}
\put(424,274){\circle*{4}}
\put(221,273){\circle*{4}}
\put(232,273){\circle*{4}}
\put(243,273){\circle*{4}}
\put(286,260){$v_{2}$}
\put(346,259){$v_{3}$}
\put(257,272){$e_{4}$}
\put(317,269){$e'$}
\put(368,272){$e_{3}$}
\put(489,166){$v_{1}$}
\put(463,178){$e_{1}$}
\put(434,166){$v_{2}$}
\put(407,177){$e_{4}$}
\put(516,178){$e_{2}$}
\put(538,165){$v_{3}$}
\put(563,176){$e_{3}$}
\put(460,43){$v_{1}$}
\put(433,57){$e_{1}$}
\put(407,43){$v_{2}$}
\put(379,56){$e_{4}$}
\put(489,54){$e'$}
\put(514,46){$u$}
\put(535,55){$e'_{2}$}
\put(562,45){$v_{3}$}
\put(590,56){$e_{3}$}
\put(342,58){\circle*{4}}
\put(353,58){\circle*{4}}
\put(364,58){\circle*{4}}
\put(593,177){\circle*{4}}
\put(604,177){\circle*{4}}
\put(615,177){\circle*{4}}
\put(368,177){\circle*{4}}
\put(379,177){\circle*{4}}
\put(390,177){\circle*{4}}
\put(619,57){\circle*{4}}
\put(630,57){\circle*{4}}
\put(641,57){\circle*{4}}
\qbezier(321,368)(321,345)(321,321)
\qbezier(332,368)(332,345)(332,321)
\qbezier(313,334)(320,323)(327,312)
\qbezier(327,312)(333,324)(340,336)
\qbezier(491,147)(491,124)(491,101)
\qbezier(503,147)(503,124)(503,101)
\qbezier(482,113)(490,102)(498,91)
\qbezier(498,91)(504,103)(511,114)
\put(317,230){$\mathcal{G}_{0}$}
\put(494,17){$\mathcal{G}_{2}$}
\put(89,180){\circle*{4}}
\put(153,180){\circle*{4}}
\put(213,179){\circle*{4}}
\qbezier(583,115)(584,115)(583,115)
\qbezier(145,146)(145,123)(145,99)
\qbezier(156,146)(156,123)(156,99)
\qbezier(137,112)(144,101)(151,90)
\qbezier(165,112)(158,101)(151,90)
\put(72,58){\circle*{4}}
\put(126,58){\circle*{4}}
\put(178,58){\circle*{4}}
\put(230,58){\circle*{4}}
\put(274,57){\circle*{4}}
\put(285,57){\circle*{4}}
\put(296,57){\circle*{4}}
\put(60,179){$e_{4}$}
\put(86,168){$v_{2}$}
\put(118,178){$e_{1}$}
\put(144,168){$v_{1}$}
\put(180,179){$e_{2}$}
\put(207,167){$v_{3}$}
\put(227,180){$e_{3}$}
\put(42,57){$e_{4}$}
\put(67,45){$v_{2}$}
\put(93,55){$e'_{1}$}
\put(123,46){$u$}
\put(148,54){$e'$}
\put(171,45){$v_{1}$}
\put(199,56){$e_{2}$}
\put(224,45){$v_{3}$}
\put(0,58){\circle*{4}}
\put(11,58){\circle*{4}}
\put(22,58){\circle*{4}}
\put(24,180){\circle*{4}}
\put(35,180){\circle*{4}}
\put(46,180){\circle*{4}}
\put(248,57){$e_{3}$}
\put(253,180){\circle*{4}}
\put(264,180){\circle*{4}}
\put(275,180){\circle*{4}}
\put(138,15){$\mathcal{G}_{1}$}
\qbezier(267,404)(301,429)(331,404)
\qbezier(267,404)(300,381)(326,400)
\qbezier(331,404)(358,431)(387,405)
\qbezier(331,404)(364,383)(387,405)
\qbezier(387,405)(391,419)(425,421)
\qbezier(387,405)(392,388)(426,388)
\qbezier(228,420)(266,419)(267,404)
\qbezier(267,404)(265,388)(227,387)
\qbezier(290,273)(281,287)(253,290)
\qbezier(252,257)(279,259)(290,273)
\qbezier(290,273)(319,296)(352,273)
\put(351,273){\circle*{4}}
\qbezier(290,273)(325,252)(351,273)
\qbezier(351,273)(370,255)(388,259)
\qbezier(351,273)(362,291)(388,292)
\qbezier(89,180)(124,205)(153,180)
\qbezier(89,180)(127,156)(153,180)
\qbezier(153,180)(186,206)(213,179)
\qbezier(153,180)(185,160)(213,179)
\qbezier(59,197)(82,194)(89,180)
\qbezier(58,166)(79,166)(89,180)
\qbezier(213,179)(218,193)(247,196)
\qbezier(213,179)(225,162)(249,166)
\qbezier(437,178)(469,203)(495,178)
\qbezier(437,178)(469,160)(495,178)
\qbezier(495,178)(524,200)(548,177)
\qbezier(495,178)(522,160)(548,177)
\qbezier(401,164)(425,163)(437,178)
\qbezier(437,178)(432,193)(403,194)
\qbezier(548,177)(554,191)(581,195)
\qbezier(548,177)(566,160)(582,164)
\qbezier(72,58)(103,80)(126,58)
\qbezier(72,58)(102,37)(126,58)
\qbezier(126,58)(151,80)(178,58)
\qbezier(126,58)(152,39)(178,58)
\qbezier(178,58)(205,77)(230,58)
\qbezier(178,58)(201,35)(230,58)
\qbezier(32,73)(66,73)(72,58)
\qbezier(32,45)(57,41)(72,58)
\qbezier(230,58)(244,73)(267,73)
\qbezier(230,58)(249,42)(268,44)
\put(465,56){\circle*{4}}
\qbezier(411,57)(440,80)(465,56)
\qbezier(411,57)(442,39)(465,56)
\put(516,57){\circle*{4}}
\qbezier(465,56)(494,80)(516,57)
\qbezier(465,56)(495,41)(516,57)
\put(568,57){\circle*{4}}
\qbezier(516,57)(545,79)(568,57)
\qbezier(516,57)(546,40)(568,57)
\qbezier(376,74)(402,74)(411,57)
\qbezier(411,57)(397,42)(376,43)
\qbezier(568,57)(575,71)(604,73)
\qbezier(568,57)(584,41)(606,43)
\end{picture}
\end{center}

\begin{lemma} \label{le3,13}
Suppose $\mathcal{G}$ is a connected hypergraph with spectral radius $\rho$ and principal eigenvector $X$. $e_{1}$, $e_{2}$, $e_{3}$, $e_{4}$ are edges in $\mathcal{G}$, where $|e_{1}|, |e_{2}|, |e_{3}|, |e_{4}|\geq 3$, $e_{1}\cap e_{2}=\{v_{1}\}$, $e_{1}\cap e_{4}=\{v_{2}\}$, $e_{2}\cap e_{3}=\{v_{3}\}$, $deg_{\mathcal{G}}(v_{1})=deg_{\mathcal{G}}(v_{2})=deg_{\mathcal{G}}(v_{3})=2$,
$deg_{\mathcal{G}}(v)=1$ for $v\in (e_{1}\cup e_{2})\setminus \{v_{1}, v_{2}, v_{3}\}$.

(1) Let $e^{'}\subset (e_{1}\cup e_{2})$ satisfy that $\{v_{2}, v_{3}\}\subseteq e^{'}$, $e^{'}\notin E(\mathcal{G})$. Let $\mathcal{G}_{0}=\mathcal{G}-e_{1}-e_{2}+e^{'}$ and $t=|e^{'}|$ (see Fig. 3.2).

(1.1) If $t\geq \max\{|e_{1}|, |e_{2}|\}$, $x_{v_{1}}\geq x_{v_{2}}$, $x_{v_{1}}\geq x_{v_{3}}$, then $\rho(\mathcal{G}_{0})\leq\rho(\mathcal{G})$ with equality if and only if $|e^{'}|=|e_{1}|= |e_{2}|$ and $x_{v_{1}}=x_{v_{2}}=x_{v_{3}}$. Moreover, if $t> \max\{|e_{1}|, |e_{2}|\}$, $x_{v_{1}}\geq x_{v_{2}}$, $x_{v_{1}}\geq x_{v_{3}}$, then $\rho(\mathcal{G}_{0})<\rho(\mathcal{G})$.

(1.2) If $t\leq \max\{|e_{1}|, |e_{2}|\}$, $x_{v_{1}}\leq x_{v_{2}}$, $x_{v_{1}}\leq x_{v_{3}}$, then $\rho(\mathcal{G}_{0})\geq\rho(\mathcal{G})$ with equality if and only if $|e^{'}|=|e_{1}|= |e_{2}|$ and $x_{v_{1}}=x_{v_{2}}=x_{v_{3}}$. Moreover, if $t< \max\{|e_{1}|, |e_{2}|\}$, $x_{v_{1}}\leq x_{v_{2}}$, $x_{v_{1}}\leq x_{v_{3}}$, then $\rho(\mathcal{G}_{0})>\rho(\mathcal{G})$.

(2) Let $e^{'}_{1}=(e_{1}\setminus \{v_{1}\})\cup \{u\}$, $e^{'}_{2}=(e_{2}\setminus \{v_{1}\})\cup \{u\}$, $e^{'}=\{v_{1}$, $u_{1}$, $u_{2}$, $\ldots$, $u_{t-2}$, $u\}$ where $u\notin V(\mathcal{G})$, $u_{i}\notin V(\mathcal{G})$ for $1\leq i\leq t-2$, $\mathcal{G}_{1}=\mathcal{G}-e_{1}+e^{'}_{1}+e^{'}$, $\mathcal{G}_{2}=\mathcal{G}-e_{2}+e^{'}_{2}+e^{'}$ (see Fig. 3.2).

(2.1) If $t\leq \min\{|e_{1}|$, $|e_{2}|\}$, $x_{v_{1}}\geq x_{v_{2}}$, $x_{v_{1}}\geq x_{v_{3}}$, then $\rho(\mathcal{G}_{1})\geq\rho(\mathcal{G})$, $\rho(\mathcal{G}_{2})\geq\rho(\mathcal{G})$ with either equality holding if and only if $|e^{'}|=|e_{1}|= |e_{2}|$, $x_{v_{1}}=x_{v_{2}}=x_{v_{3}}=x_{u}$ and $x_{z}=x_{\omega}$ for $z,\omega\in (e_{1}\cup e_{2}\cup e^{'})\setminus\{v_{1}, v_{2}, u\}$. Moreover, if $t< \min\{|e_{1}|$, $|e_{2}|\}$, $x_{v_{1}}\geq x_{v_{2}}$, $x_{v_{1}}\geq x_{v_{3}}$, then $\rho(\mathcal{G}_{1})>\rho(\mathcal{G})$, $\rho(\mathcal{G}_{2})>\rho(\mathcal{G})$.

(2.2) If $t\geq \max\{e_{1}$, $e_{2}\}$, $x_{v_{1}}\leq x_{v_{2}}$, $x_{v_{1}}\leq x_{v_{3}}$, then $\rho(\mathcal{G}_{1})\leq\rho(\mathcal{G})$, $\rho(\mathcal{G}_{2})\leq\rho(\mathcal{G})$ with either equality holding if and only if $|e^{'}|=|e_{1}|= |e_{2}|$, $x_{v_{1}}=x_{v_{2}}=x_{v_{3}}=x_{u}$ and $x_{z}=x_{\omega}$ for $z,\omega\in (e_{1}\cup e_{2}\cup e^{'})\setminus\{v_{1}, v_{2}, u\}$. Moreover, if $t> \max\{e_{1}$, $e_{2}\}$, $x_{v_{1}}\leq x_{v_{2}}$, $x_{v_{1}}\leq x_{v_{3}}$, then $\rho(\mathcal{G}_{1})<\rho(\mathcal{G})$, $\rho(\mathcal{G}_{2})<\rho(\mathcal{G})$.

\end{lemma}

\begin{proof}
(1.1) Suppose $e_{1}=\{v_{1}$, $v_{\alpha(1,1)}$, $v_{\alpha(1,2)}$, $\ldots$, $v_{\alpha(1,j_{1}-2)}$, $v_{2}\}$, $e_{2}=\{v_{1}$, $v_{\alpha(2,1)}$, $v_{\alpha(2,2)}$, $\ldots$, $v_{\alpha(2,j_{2}-2)}$, $v_{3}\}$. By Lemma \ref{le3,08}, we have $x_{v_{\alpha(1,w)}}=x_{v_{\alpha(1,z)}}< \min\{x_{v_{1}}, x_{v_{2}}\}$ for $1\leq w< z\leq j_{1}-2$, $x_{v_{\alpha(2,w)}}=x_{v_{\alpha(2,z)}}< \min\{x_{v_{1}}, x_{v_{3}}\}$ for $1\leq w< z\leq j_{2}-2$. Let $Y$ be a vector on $\mathcal{G}_{0}$ satisfying that
$$\left \{\begin{array}{ll}
 y_{v}=\min \{x_{z}: z\in (e_{1}\cup e_{2})\setminus \{v_{1}, v_{2}, v_{3}\}\},\ & \ v\in e^{'}\setminus\{v_{2}, v_{3}\}\\
  y_{v}=x_{v},\ & \ others. \end{array}\right.$$

Note that $|e^{'}|\geq \max\{|e_{1}|, |e_{2}|\}$, $x(v_{1})\geq x(v_{2})$, $x(v_{1})\geq x(v_{3})$. Without loss of generality, suppose $\min \{x_{v}: v\in (e_{1}\cup e_{2})\}=x_{v_{\alpha(2,1)}}$.
For $v\in (e^{'}\setminus \{v_{2}, v_{3}\})$, noting that $deg_{\mathcal{G}_{0}}(v)=1$ and $x_{v_{\alpha(2,1)}}< \min\{x_{v_{1}}, x_{v_{3}}\}$, we have $$(\mathcal {A}_{\mathcal{G}_{0}}Y)_{v}=\frac{(t-3)y_{v}+y_{v_{2}}+y_{v_{3}}}{t-1}=\frac{(t-3)x_{v_{\alpha(2,1)}}+x_{v_{2}}+x_{v_{3}}}{t-1}\hspace{3cm}$$
$$=\frac{(j_{2}-3)x_{v_{\alpha(2,1)}}+x_{v_{2}}+x_{v_{3}}+(t-j_{2})x_{v_{\alpha(2,1)}}}{j_{2}-1+t-j_{2}}\hspace{2.4cm}$$
$$\leq\frac{(j_{2}-3)x_{v_{\alpha(2,1)}}+x_{v_{1}}+x_{v_{3}}+(t-j_{2})x_{v_{\alpha(2,1)}}}{j_{2}-1+t-j_{2}}\hspace{2.4cm}$$
$$\leq\frac{(j_{2}-3)x_{v_{\alpha(2,1)}}+x_{v_{1}}+x_{v_{3}}}{j_{2}-1}\hspace{2.5cm}(by\ Lemma\ \ref{le3,07,01})$$
$$= \rho x_{v_{\alpha(2,1)}}=\rho y_{v}.\hspace{6.8cm}$$ In the same way, we get
$$(\mathcal {A}_{\mathcal{G}_{0}}Y)_{v_{2}}=\frac{(t-2)x_{v_{\alpha(2,1)}}+y_{v_{3}}}{t-1}=\frac{(t-2)x_{v_{\alpha(2,1)}}+x_{v_{3}}}{t-1}\leq \rho x_{v_{2}}=\rho y_{v_{2}};$$
$$(\mathcal {A}_{\mathcal{G}_{0}}Y)_{v_{3}}=\frac{(t-2)x_{v_{\alpha(2,1)}}+y_{v_{2}}}{t-1}=\frac{(t-2)x_{v_{\alpha(2,1)}}+x_{v_{2}}}{t-1}\leq \rho x_{v_{3}}=\rho y_{v_{3}};$$
for $v\in (V(\mathcal{G}_{0})\setminus e^{'})$, $(\mathcal {A}_{\mathcal{G}_{0}}Y)_{v}=(\mathcal {A}_{\mathcal{G}}X)_{v}=\rho y_{v}$.
By lemma \ref{le3,04}, it follows that $\rho(\mathcal{G}_{0})\leq\rho(\mathcal{G})$. Note that $Y$ is positive. Combining Lemma \ref{le3,04,01}, we find that if $\rho(\mathcal{G}_{0})=\rho(\mathcal{G})$, then $\mathcal {A}_{\mathcal{G}_{0}}Y=\rho Y$. Thus it follows that $|e^{'}|=|e_{1}|= |e_{2}|$ and $x_{v_{1}}=x_{v_{2}}=x_{v_{3}}$. Conversely, if $|e^{'}|=|e_{1}|= |e_{2}|$ and $x_{v_{1}}=x_{v_{2}}=x_{v_{3}}$, then it can be checked as above that for $v\in (e^{'}\setminus \{v_{2}, v_{3}\})$, $(\mathcal {A}_{\mathcal{G}_{0}}Y)_{v}=\rho y_{v}$; $(\mathcal {A}_{\mathcal{G}_{0}}Y)_{v_{2}}=\rho y_{v_{2}};$
$(\mathcal {A}_{\mathcal{G}_{0}}Y)_{v_{3}}=\rho y_{v_{3}};$
for $v\in (V(\mathcal{G}_{0})\setminus e^{'})$, $(\mathcal {A}_{\mathcal{G}_{0}}Y)_{v}=(\mathcal {A}_{\mathcal{G}}X)_{v}=\rho y_{v}$. Thus it follows that $\rho(\mathcal{G}_{0})=\rho(\mathcal{G})$.

If $t> \max\{|e_{1}|, |e_{2}|\}$, $x_{v_{1}}\geq x_{v_{2}}$, $x_{v_{1}}\geq x_{v_{3}}$, combining Lemma \ref{le3,07,01}, as the above proof, we get $\rho(\mathcal{G}_{0})<\rho(\mathcal{G})$.

From the above narrations, then (1.1) follows as desired.

(1.2) Let $Y$ be a vector on $\mathcal{G}_{0}$ satisfying that
$$\left \{\begin{array}{ll}
 y(v)=\max \{x_{z}: z\in (e_{1}\cup e_{2})\setminus \{v_{1}, v_{2}, v_{3}\}\},\ & \ v\in e^{'}\setminus\{v_{2}, v_{3}\}\\
  y(v)=x_{v},\ & \ others. \end{array}\right.$$
Then (1.2) is proved as (1.1).

(2.1) For both $\mathcal{G}_{1}$ and $\mathcal{G}_{2}$, let $Y$ be a vector satisfying that
$$\left \{\begin{array}{ll}
y(u)=x_{v_{1}}\ & \ \\
 y(v)=\max \{x_{z}: z\in (e_{1}\cup e_{2})\setminus \{v_{1}, v_{2}, v_{3}\}\},\ & \ v\in e^{'}\setminus \{v_{1}, u\}\\
  y(v)=x_{v},\ & \ others. \end{array}\right.$$
Then (2.1) is proved as (1.1).

(2.2) For both $\mathcal{G}_{1}$ and $\mathcal{G}_{2}$, let $Y$ be a vector satisfying that
$$\left \{\begin{array}{ll}
y(u)=x_{v_{1}}\ & \ \\
 y(v)=\min \{x_{z}: z\in (e_{1}\cup e_{2})\setminus \{v_{1}, v_{2}, v_{3}\}\},\ & \ v\in e^{'}\setminus \{v_{1}, u\}\\
  y(v)=x_{v},\ & \ others. \end{array}\right.$$
Then (2.2) is proved as (1.1).

 Thus the result follows. This completes the proof. \ \ \ \ \ $\Box$
\end{proof}

\begin{lemma} \label{le3,14}
Let $e$ be a new edge not containing in connected hypergraph $\mathcal{G}$. Let $\mathcal{G}^{'}=\mathcal{G}+e$. If $\mathcal{G}^{'}$ is also connected, then $\rho(\mathcal{G}^{'})>\rho(\mathcal{G})$.
\end{lemma}

\begin{proof}
Let $X$ be
the $principal$ $eigenvector$ of $\mathcal{G}$, $Y$ be a vector on $\mathcal{G}^{'}$ satisfying that
$$\left \{\begin{array}{ll}
 y_{v}=x_{v},\ & \ v\in V(\mathcal{G})\\
 \\
  y_{v}=0,\ & \ others. \end{array}\right.$$ Then $Y^{T}\mathcal {A}_{\mathcal{G}^{'}}Y-X^{T}\mathcal {A}_{\mathcal{G}}X\geq 0$, $Y^{T}Y=X^{T}X$. It follows that $\rho(\mathcal{G}^{'})\geq\rho(\mathcal{G})$.
  Suppose that $\rho(\mathcal{G}^{'})=\rho(\mathcal{G})$. Then $\rho(\mathcal{G}^{'})=Y^{T}\mathcal {A}_{\mathcal{G}^{'}}Y=X^{T}\mathcal {A}_{\mathcal{G}}X=\rho(\mathcal{G})$, and then $Y$ is a principal eigenvector of $\mathcal{G}^{'}$. If there exists $y_{v}=0$, then we get a contradiction because the principal eigenvector of $\mathcal{G}^{'}$ is positive by Perron-Frobenius theorem.

Suppose $Y$ is positive next. Note that $Y=X$ if $Y$ is positive. It follows that $V(\mathcal{G}^{'})=V(\mathcal{G})$ now. Denote by $e=\{v_{1}$, $v_{2}$, $\ldots$, $v_{k}\}$. Then $e\subseteq V(\mathcal{G})$, and
$$\rho(\mathcal{G}^{'})y_{v_{1}}=(\mathcal {A}_{\mathcal{G}^{'}}Y)_{v_{1}}=(\mathcal {A}_{\mathcal{G}^{'}}X)_{v_{1}}+\frac{1}{k-1}\sum_{i=2}^{k} x_{v_{i}}=\rho(\mathcal{G})x_{v_{1}}+\frac{1}{k-1}\sum_{i=2}^{k} x_{v_{i}}> \rho(\mathcal{G})x_{v_{1}},$$ which contradicts $\rho(\mathcal{G}^{'})=\rho(\mathcal{G})$. As a result, we get that $\rho(\mathcal{G}^{'})>\rho(\mathcal{G})$.
This completes the proof. \ \ \ \ \ $\Box$
\end{proof}

Denote by $\mathcal{G}(\mathcal{D}v; p, q; v_{p+q}\mathcal{H})$ the $k$-uniform connected hypergraph obtained from $k$-uniform hypergraph $\mathcal{D}$ and $k$-uniform hypergraph $\mathcal{H}$ by adding a pendant path $P_{1}$ with length $p$ at vertex $v$ of $\mathcal{D}$, and adding a path $P_{2}$ with length $q$ between vertex $v$ and vertex $v_{p+q}$ of $\mathcal{H}$, where $\mathcal{D}$ and $\mathcal{H}$ are two disjoint, $V(P_{1})\cap V(\mathcal{D})=\{v\}$, $V(P_{2})\cap V(\mathcal{D})=\{v\}$, $V(P_{2})\cap V(\mathcal{H})=\{v_{p+q}\}$ (see two examples in Fig. 3.3). In particular, if $H=v_{p+q}$, we denote by $\mathcal{G}(\mathcal{D}v; p, q; v_{p+q})$ for $\mathcal{G}(\mathcal{D}v; p, q; v_{p+q}H)$ for short.

\

\setlength{\unitlength}{0.6pt}
\begin{center}
\begin{picture}(765,153)
\qbezier(0,97)(0,108)(9,116)\qbezier(9,116)(19,124)(33,124)\qbezier(33,124)(46,124)(56,116)\qbezier(56,116)(66,108)(66,97)
\qbezier(66,97)(66,85)(56,77)\qbezier(56,77)(46,70)(33,70)\qbezier(33,70)(19,70)(9,77)\qbezier(9,77)(0,85)(0,97)
\put(66,98){\circle*{4}}
\put(110,130){\circle*{4}}
\put(133,130){\circle*{4}}
\put(143,130){\circle*{4}}
\put(123,130){\circle*{4}}
\put(155,130){\circle*{4}}
\put(210,130){\circle*{4}}
\qbezier(211,131)(211,140)(219,146)\qbezier(219,146)(227,153)(239,153)\qbezier(239,153)(250,153)(258,146)\qbezier(258,146)(267,140)(267,131)
\qbezier(267,131)(267,121)(258,115)\qbezier(258,115)(250,109)(239,109)\qbezier(239,109)(227,109)(219,115)\qbezier(219,115)(211,121)(211,131)
\put(110,61){\circle*{4}}
\put(184,61){\circle*{4}}
\put(194,61){\circle*{4}}
\put(174,61){\circle*{4}}
\put(204,61){\circle*{4}}
\put(257,61){\circle*{4}}
\put(18,92){$\mathcal{D}$}
\put(55,94){$v$}
\put(107,68){$v_{1}$}
\put(260,58){$v_{p}$}
\put(95,138){$v_{p+1}$}
\put(214,126){$v_{p+q}$}
\put(250,127){$\mathcal{H}$}
\put(71,20){$\mathcal{G}(Dv; p, q; v_{p+q}\mathcal{H})$}
\qbezier(306,96)(306,107)(315,115)\qbezier(315,115)(325,123)(339,123)\qbezier(339,123)(352,123)(362,115)\qbezier(362,115)(372,107)(372,96)
\qbezier(372,96)(372,84)(362,76)\qbezier(362,76)(352,69)(339,69)\qbezier(339,69)(325,69)(315,76)\qbezier(315,76)(306,84)(306,96)
\put(372,95){\circle*{4}}
\put(420,95){\circle*{4}}
\put(479,95){\circle*{4}}
\put(491,95){\circle*{4}}
\put(502,95){\circle*{4}}
\put(513,96){\circle*{4}}
\put(561,96){\circle*{4}}
\put(609,96){\circle*{4}}
\put(621,96){\circle*{4}}
\put(633,96){\circle*{4}}
\put(644,96){\circle*{4}}
\put(657,96){\circle*{4}}
\put(709,96){\circle*{4}}
\qbezier(709,100)(709,109)(717,116)\qbezier(717,116)(725,124)(737,124)\qbezier(737,124)(748,124)(756,116)\qbezier(756,116)(765,109)(765,100)
\qbezier(765,100)(765,90)(756,83)\qbezier(756,83)(748,76)(737,76)\qbezier(737,76)(725,76)(717,83)\qbezier(717,83)(709,90)(709,100)
\put(471,21){$\mathcal{G}(\mathcal{D}v; 0, p+q; v_{p+q}\mathcal{H})$}
\put(205,-9){Fig. 3.3. $\mathcal{G}(Dv; p, q; v_{p+q}\mathcal{H})$ and $\mathcal{G}(\mathcal{D}v; 0, p+q; v_{p+q}\mathcal{H})$}
\put(329,92){$\mathcal{D}$}
\put(746,96){$\mathcal{H}$}
\put(361,92){$v$}
\put(553,105){$v_{p}$}
\put(600,104){$v_{p+1}$}
\put(712,93){$v_{p+q}$}
\put(413,102){$v_{1}$}
\put(387,80){$e_{1}$}
\put(466,95){\circle*{4}}
\put(436,78){$e_{2}$}
\put(459,102){$v_{2}$}
\put(531,80){$e_{p}$}
\put(162,61){\circle*{4}}
\put(156,66){$v_{2}$}
\put(68,66){$e_{1}$}
\put(126,45){$e_{2}$}
\put(224,45){$e_{p}$}
\put(55,125){$e_{p+1}$}
\put(572,82){$e_{p+1}$}
\qbezier(66,98)(81,126)(110,130)
\qbezier(110,130)(99,105)(66,98)
\qbezier(155,130)(181,143)(210,130)
\qbezier(155,130)(180,119)(210,130)
\qbezier(66,98)(91,89)(110,61)
\qbezier(66,98)(76,69)(110,61)
\qbezier(110,61)(135,73)(162,61)
\qbezier(110,61)(135,52)(162,61)
\qbezier(204,61)(231,72)(257,61)
\qbezier(204,61)(231,52)(257,61)
\qbezier(372,95)(396,105)(420,95)
\qbezier(372,95)(396,86)(420,95)
\qbezier(420,95)(444,106)(466,95)
\qbezier(420,95)(444,87)(466,95)
\qbezier(513,96)(538,108)(561,96)
\qbezier(513,96)(539,89)(561,96)
\qbezier(561,96)(585,109)(609,96)
\qbezier(561,96)(586,88)(609,96)
\qbezier(657,96)(683,109)(709,96)
\qbezier(657,96)(685,89)(709,96)
\end{picture}
\end{center}

\begin{lemma} \label{le3,39}
If $p, q> 0$, then
$\rho(\mathcal{G}(\mathcal{D}v; p, q; v_{p+q}\mathcal{H})) > \rho(\mathcal{G}(\mathcal{D}v; 0, p+q; v_{p+q}\mathcal{H}))$.
\end{lemma}

\begin{proof}
Let $X$ be
the principal eigenvector of the uniform hypergraph $\mathcal{G}((\mathcal{D}v; 0, p+q; v_{p+q}\mathcal{H}))$. Assume that in $\mathcal{D}$, the edges incident with $v$ are $\varepsilon_{1}$, $\varepsilon_{2}$, $\ldots$, $\varepsilon_{\eta}$. Let $e^{'}_{p+1}=(e_{p+1}\setminus\{v_{p}\})\cup \{v\}$, $\mathcal{G}_{1}=\mathcal{G}((\mathcal{D}v; 0, p+q; v_{p+q}\mathcal{H}))-e_{p+1}+e^{'}_{p+1}$. If $x_{v}\geq x_{v_{p}}$, then $\rho(\mathcal{G}_{1}) > \rho(\mathcal{G}(\mathcal{D}v; 0, p+q; v_{p+q}\mathcal{H}))$ by Lemma \ref{le03.02}. Let $\varepsilon^{'}_{i}=(\varepsilon_{i}\setminus \{v\})\cup \{v_{p}\}$ for $1\leq i\leq \eta$, $\mathcal{G}_{2}=\mathcal{G}((\mathcal{D}v; 0, p+q; v_{p+q}\mathcal{H}))-\sum^{\eta}_{i=1} \varepsilon_{i}+\sum^{\eta}_{i=1} \varepsilon^{'}_{i}$. If $x_{v}\leq x_{v_{p}}$, then $\rho(\mathcal{G}_{2}) > \rho(\mathcal{G}(\mathcal{D}v; 0, p+q; v_{p+q}\mathcal{H}))$ by Lemma \ref{le03.02}. Note that both $G_{1}$ and $G_{2}$ are isomorphic
to $\mathcal{G}(\mathcal{D}v; p, q; v_{p+q}\mathcal{H})$. Thus we get that $\rho(\mathcal{G}(\mathcal{D}v; p, q; v_{p+q}\mathcal{H})) > \rho(\mathcal{G}(\mathcal{D}v; 0, p+q; v_{p+q}\mathcal{H}))$.
This completes the proof. \ \ \ \ \ $\Box$
\end{proof}

\setlength{\unitlength}{0.65pt}
\begin{center}
\begin{picture}(691,114)
\qbezier(55,45)(55,56)(64,64)\qbezier(64,64)(74,72)(88,72)\qbezier(88,72)(101,72)(111,64)\qbezier(111,64)(121,56)(121,45)\qbezier(121,45)(121,33)(111,25)
\qbezier(111,25)(101,18)(88,18)\qbezier(88,18)(74,18)(64,25)\qbezier(64,25)(55,33)(55,45)
\put(121,45){\circle*{4}}
\put(172,45){\circle*{4}}
\put(184,44){\circle*{4}}
\put(196,44){\circle*{4}}
\put(207,44){\circle*{4}}
\put(377,46){\circle*{4}}
\put(432,46){\circle*{4}}
\put(489,46){\circle*{4}}
\put(500,45){\circle*{4}}
\put(512,45){\circle*{4}}
\put(523,45){\circle*{4}}
\put(534,46){\circle*{4}}
\put(590,46){\circle*{4}}
\put(260,-9){Fig. 3.4. $\mathcal{G}(\mathcal{D}v_{0}; 0, t; v_{t})$}
\put(75,40){$\mathcal{D}$}
\put(104,43){$v_{0}$}
\put(426,54){$v_{p}$}
\put(400,30){$e_{p}$}
\put(450,30){$e_{p+1}$}
\put(476,54){$v_{p+1}$}
\put(594,45){$v_{t}$}
\put(167,52){$v_{1}$}
\put(220,45){\circle*{4}}
\put(276,45){\circle*{4}}
\put(330,45){\circle*{4}}
\put(364,44){\circle*{4}}
\put(353,44){\circle*{4}}
\put(341,44){\circle*{4}}
\put(269,54){$v_{\zeta}$}
\put(317,53){$v_{\zeta+1}$}
\put(210,55){$v_{\zeta-1}$}
\put(365,55){$v_{p-1}$}
\put(518,54){$v_{t-1}$}
\put(143,27){$e_{1}$}
\put(235,29){$e_{\zeta-1}$}
\put(288,29){$e_{\zeta+1}$}
\put(558,31){$e_{t}$}
\qbezier(121,45)(147,57)(172,45)
\qbezier(121,45)(146,34)(172,45)
\qbezier(220,45)(250,58)(276,45)
\qbezier(220,45)(250,34)(276,45)
\qbezier(276,45)(303,58)(330,45)
\qbezier(276,45)(303,36)(330,45)
\qbezier(377,46)(407,58)(432,46)
\qbezier(432,46)(460,58)(489,46)
\qbezier(534,46)(561,56)(590,46)
\qbezier(377,46)(408,34)(432,46)
\qbezier(432,46)(459,34)(489,46)
\qbezier(534,46)(562,35)(590,46)
\end{picture}
\end{center}

\begin{lemma} \label{le3,40}
Let $X$ be
the $principal$ $eigenvector$ of the uniform hypergraph $\mathcal{G}(\mathcal{D}v_{0}; 0, t; v_{t})$. Denote by $\mathcal{P}=v_{0}e_{1}v_{1}e_{2}v_{2}\cdots e_{t}v_{t}$ the pendant path from vertex $v_{0}$ to vertex $v_{t}$ in $\mathcal{G}(\mathcal{D}v_{0}; 0, t; v_{t})$ and $e_{i}=\{v_{i-1}$, $v_{a(i,1)}$, $v_{a(i,2)}$, $\ldots$, $v_{a(i,k-2)}$, $v_{i}\}$ for $1\leq i\leq t$ (see Fig. 3.4). Suppose $x_{v_{p}}=\max\{x_{v_{i}}: 0\leq i\leq t\}$. Then

$\mathrm{(1)}$ $p\leq t-1$.

Let $t=p+q$. Moreover, we have

$\mathrm{(2)}$ if $p> 0$, then $x_{v_{i}}\leq x_{v_{i+1}}$ for $0\leq i\leq p-1$, $x_{v_{i}}\geq x_{v_{i+1}}$ for $p\leq i\leq t-1$; if $p= 0$, then $x_{v_{i}}\geq x_{v_{i+1}}$ for $0\leq i\leq t-1$.

$\mathrm{(3)}$ if $p> 0$ and there exists $\omega\leq p$ and $\eta\leq q$ such that $x_{v_{p-\omega}}\geq x_{v_{p+\eta}}$, then

$\mathrm{(3.1)}$ if $\omega\leq \eta$, then $x_{v_{p-\omega+i}}\geq x_{v_{p+\eta-i}}$ for $0\leq i\leq \omega$, $x_{v_{a(p-\omega+i,1)}}\geq x_{v_{a(p+\eta-i+1,1)}}$ for $1\leq i\leq \omega$.

$\mathrm{(3.2)}$ if $\omega\geq \eta$, then $x_{v_{p-\omega+i}}\geq x_{v_{p+\eta-i}}$ for $0\leq i\leq \eta-1$, $x_{v_{j}}=x_{v_{p}}$ for $p-\omega+\eta\leq j\leq p$, and $x_{v_{a(p-\omega+i,1)}}\geq x_{v_{a(p+\eta-i+1,1)}}$ for $1\leq i\leq \eta$.

$\mathrm{(4)}$ if $p> 0$ and there exists $\omega\leq p$ and $\eta\leq q$ such that $x_{v_{p-\omega}}\leq x_{v_{p+\eta}}$, then

$\mathrm{(4.1)}$ if $\omega\leq \eta$, then $x_{v_{p-\omega+i}}\leq x_{v_{p+\eta-i}}$ for $0\leq i\leq \omega-1$, $x_{v_{j}}=x_{v_{p}}$ for $p+1\leq j\leq p+\eta-\omega$, and $x_{v_{a(p-\omega+i,1)}}\leq x_{v_{a(p+\eta-i+1,1)}}$ for $1\leq i\leq \omega$.

$\mathrm{(4.2)}$ if $\omega\geq \eta$, then $x_{v_{p-\omega+i}}\leq x_{v_{p+\eta-i}}$ for $0\leq i\leq \eta$, $x_{v_{a(p-\omega+i,1)}}\leq x_{v_{a(p+\eta-i+1,1)}}$ for $1\leq i\leq \eta$.

$\mathrm{(5)}$ if $p> 0$ and there exists $\omega\leq p$, $\eta\leq q$ such that $x_{v_{p-\omega}}=x_{v_{p+\eta}}$, then

$\mathrm{(5.1)}$ if $\omega\leq \eta$, then $x_{v_{p-\omega+i}}= x_{v_{p+\eta-i}}$ for $0\leq i\leq \omega-1$, $x_{v_{j}}=x_{v_{p}}$ for $p+1\leq j\leq p+\eta-\omega$, and $x_{v_{a(p-\omega+i,1)}}= x_{v_{a(p+\eta-i+1,1)}}$ for $1\leq i\leq \omega$.

$\mathrm{(5.2)}$ if $\omega\geq \eta$, then $x_{v_{p-\omega+i}}= x_{v_{p+\eta-i}}$ for $0\leq i\leq \eta-1$, $x_{v_{j}}=x_{v_{p}}$ for $p-\omega+\eta\leq j\leq p$, $x_{v_{a(p-\omega+i,1)}}= x_{v_{a(p+\eta-i+1,1)}}$ for $1\leq i\leq \eta$.

\end{lemma}

\begin{proof}
(1) By Lemma \ref{le3,08,02}, it follows that $x_{v_{t}}<x_{v_{t-1}}$. Thus $p\leq t-1$.

(2) Using Lemma \ref{le3,08} gets that for $1\leq i\leq t-1$, $$x_{v_{a(i,1)}}=x_{v_{a(i,2)}}= \cdots = x_{v_{a(i,k-2)}}< \min\{x_{v_{i-1}}, x_{v_{i}}\}.$$

Note that $\mathcal{G}(\mathcal{D}v_{0}; 0, t; v_{t})$ is uniform and  $p\leq t-1$.

{\bf Case 1} $p> 0$.

{\bf Claim} $x_{v_{0}}\leq x_{v_{1}}\leq x_{v_{2}}\leq\cdots \leq x_{v_{p}}$. If $p=1$, this claim hold naturally.

For $p\geq 2$, we prove this claim by contradiction. Suppose that $x_{v_{z}}$ ($0\leq z\leq p-1$) is the first vertex from $0$ to $p-1$ such that $x_{v_{z}}> x_{v_{z+1}}$. Then there exists $z+1\leq \zeta\leq p-1$ such that $x_{v_{z}}> x_{v_{z+1}}\geq \cdots \geq x_{v_{\zeta}}\leq x_{v_{\zeta+1}}$. Let $e^{'}=\{v_{\zeta}$, $u_{1}$, $u_{2}$, $\ldots$, $u_{k-2}$, $u\}$, $e^{'}_{\zeta}=(e_{\zeta}\setminus \{v_{\zeta}\})\cup \{u\}$, where $u\notin V(\mathcal{G})$, $u_{i}\notin V(\mathcal{G})$ for $1\leq i\leq t-2$, $\mathcal{G}_{1}=\mathcal{G}-e_{\zeta}+e^{'}_{\zeta}+e^{'}$, $\mathcal{G}_{2}=\mathcal{G}_{1}- \{v_{a(t,1)}$, $v_{a(t,2)}$, $\ldots$, $v_{a(t,k-2)}$, $v_{t}\}$. By Lemma \ref{le3,13} and Lemma \ref{le3,14}, we get that $\rho(\mathcal{G}_{2})< \rho(\mathcal{G}_{1})\leq \rho(\mathcal{G})$, which contradicts $\rho(\mathcal{G})=\rho(\mathcal{G}_{2})$ because $\mathcal{G}_{2}\cong \mathcal{G}$. Thus our claim holds.

In the same way, it is proved that $x_{v_{p}}\geq x_{v_{p+1}}\geq x_{v_{p+2}}\geq \cdots \geq x_{v_{t-1}}$. Combining Lemma \ref{le3,08,02}, we get that $x_{v_{p}}\geq x_{v_{p+1}}\geq x_{v_{p+2}}\geq \cdots \geq x_{v_{t-1}}>x_{v_{t}}$.

{\bf Case 2} $p= 0$. As Case 1, it is proved that $x_{v_{i}}\geq x_{v_{i+1}}$ for $0\leq i\leq t-1$.
 Thus (2) follows.

(3) If $\omega\leq \eta$, we let $Y$ be
a vector on $\mathcal{G}(\mathcal{D}v_{0}; 0, t; v_{t})$ satisfying that
$$\left \{\begin{array}{ll}
y_{v_{p-\omega+i}}=\max \{x_{v_{p-\omega+i}}, x_{v_{p+\eta-i}}\}\ ~~~& \ 0\leq i\leq \eta;\\
\\
y_{v_{a(p-\omega+i,j)}}=\max \{x_{v_{a(p-\omega+i,1)}}, x_{v_{a(p+\eta-i+1,1)}}\}\ ~~~& \ 1\leq i\leq  \omega, 1\leq j\leq k-2;\\
\\
y_{v}=x_{v}\ ~~~& \ others. \end{array}\right.$$

As the proof of Lemma \ref{le3,13}, it is proved that $\mathcal {A}(\mathcal{G}(\mathcal{D}v_{0}; 0, t; v_{t}))Y\geq \rho(\mathcal{G}(\mathcal{D}v_{0}; 0, t; v_{t}))Y$. Using Lemma \ref{le3,05} gets that $\mathcal {A}(\mathcal{G}(\mathcal{D}v_{0}; 0, t; v_{t}))Y= \rho(\mathcal{G}(\mathcal{D}v_{0}; 0, t; v_{t}))Y$. Note that $\mathcal{G}(\mathcal{D}v_{0}; 0, t; v_{t})$ is connected. Then $\mathcal {A}(\mathcal{G}(\mathcal{D}v_{0}; 0, t; v_{t}))$ is irreducible. Consequently, it follows that the dimension of the eigenspace
of the eigenvalue $\rho(\mathcal{G}(\mathcal{D}v_{0}; 0, t; v_{t}))$ is one. Then $Y=lX$ for some $l> 0$. Then (3.1) follows.

(3.2), (4.1) and (4.2) are proved in the same way.
(5) follows from (3) or (4).

This completes the proof. \ \ \ \ \ $\Box$
\end{proof}

\begin{lemma} \label{le3,41,0}
$L_{0}, L_{1}, L_{2},\ldots, L_{f}$ ($f\geq 1$) are positive integers satisfying $L_{1}\leq L_{2}\leq\cdots\leq L_{f}\leq  L_{0}-1$ and $\sum^{f}_{i=0}L_{i}=t$. For an integer $\mu> 0$, if $t-\mu\geq L_{0}$, then there exists some $1\leq j\leq f$ such that $t-\mu> L_{1}+\cdots +L_{j}$, but $L_{0}+L_{1}+\cdots +L_{j}> t-\mu$.
\end{lemma}

\begin{proof}
If $\sum^{f}_{i=1}L_{i}< t-\mu$, then $t=\sum^{f}_{i=0}L_{i}> t-\mu$ follows naturally. Thus the result holds.

Suppose $\sum^{f}_{i=1}L_{i}\geq t-\mu$. Note that $t-\mu\geq L_{0}$ and $L_{1}\leq L_{2}\leq\cdots\leq L_{f}\leq  L_{0}-1$. Then $f\geq 2$ now, and there exists some $1\leq g\leq f-1$ such that $L_{1}+\cdots +L_{g}+L_{g+1}\geq t-\mu$, but $L_{1}+\cdots +L_{g}<t-\mu$. Note that $L_{g+1}\leq L_{0}-1$. Then it follows that $L_{0}+L_{1}+\cdots +L_{g}> t-\mu$, but $t-\mu> L_{1}+\cdots +L_{g}$.
This completes the proof. \ \ \ \ \ $\Box$
\end{proof}

\

\setlength{\unitlength}{0.7pt}
\begin{center}
\begin{picture}(533,287)
\qbezier(112,76)(112,87)(121,95)\qbezier(121,95)(131,103)(145,103)\qbezier(145,103)(158,103)(168,95)\qbezier(168,95)(178,87)(178,76)
\qbezier(178,76)(178,64)(168,56)\qbezier(168,56)(158,49)(145,49)\qbezier(145,49)(131,49)(121,56)\qbezier(121,56)(112,64)(112,76)
\put(177,76){\circle*{4}}
\put(289,75){\circle*{4}}
\put(363,73){\circle*{4}}
\put(375,73){\circle*{4}}
\put(386,73){\circle*{4}}
\put(216,-9){Fig. 3.5. $\mathcal{G}_{1}$, $\mathcal{G}_{2}$}
\put(129,74){$\mathcal{D}$}
\put(161,74){$v_{0}$}
\put(285,63){$u_{1}$}
\put(400,74){\circle*{4}}
\put(463,74){\circle*{4}}
\put(526,75){\circle*{4}}
\put(451,83){$u_{s-1}$}
\put(531,72){$u_{s}$}
\put(231,54){$\tilde{e}_{1}$}
\put(422,54){$\tilde{e}_{s-1}$}
\put(495,54){$\tilde{e}_{s}$}
\qbezier(177,76)(232,98)(289,75)
\qbezier(177,76)(224,58)(289,75)
\put(202,131){\circle*{4}}
\put(202,78){\circle*{4}}
\qbezier(202,131)(187,102)(202,78)
\qbezier(202,131)(216,111)(202,78)
\put(158,155){\circle*{4}}
\qbezier(202,131)(176,133)(158,155)
\qbezier(158,155)(197,155)(202,131)
\put(250,135){\circle*{4}}
\put(231,77){\circle*{4}}
\qbezier(250,135)(227,115)(231,77)
\qbezier(250,135)(254,96)(231,77)
\put(285,181){\circle*{4}}
\qbezier(250,135)(254,167)(285,181)
\qbezier(285,181)(283,150)(250,135)
\put(392,182){\circle*{4}}
\put(340,182){\circle*{4}}
\qbezier(392,182)(366,198)(340,182)
\qbezier(392,182)(368,169)(340,182)
\put(263,77){\circle*{4}}
\put(295,115){\circle*{4}}
\qbezier(263,77)(265,105)(295,115)
\qbezier(295,115)(295,92)(263,77)
\put(344,132){\circle*{4}}
\qbezier(295,115)(315,139)(344,132)
\qbezier(295,115)(324,111)(344,132)
\put(380,132){\circle*{4}}
\put(369,132){\circle*{4}}
\put(357,132){\circle*{4}}
\put(395,132){\circle*{4}}
\put(449,132){\circle*{4}}
\qbezier(395,132)(421,148)(449,132)
\qbezier(395,132)(421,121)(449,132)
\put(323,181){\circle*{4}}
\put(312,181){\circle*{4}}
\put(300,181){\circle*{4}}
\put(145,155){\circle*{4}}
\put(134,155){\circle*{4}}
\put(122,155){\circle*{4}}
\put(57,154){\circle*{4}}
\put(109,155){\circle*{4}}
\qbezier(57,154)(85,168)(109,155)
\qbezier(57,154)(85,139)(109,155)
\qbezier(66,259)(66,268)(75,274)\qbezier(75,274)(84,281)(97,281)\qbezier(97,281)(109,281)(118,274)\qbezier(118,274)(128,268)(128,259)
\qbezier(128,259)(128,249)(118,243)\qbezier(118,243)(109,237)(97,237)\qbezier(97,237)(84,237)(75,243)\qbezier(75,243)(66,249)(66,259)
\put(127,258){\circle*{4}}
\put(186,258){\circle*{4}}
\qbezier(127,258)(158,274)(186,258)
\qbezier(127,258)(157,247)(186,258)
\put(237,258){\circle*{4}}
\qbezier(186,258)(210,273)(237,258)
\qbezier(186,258)(213,248)(237,258)
\put(276,257){\circle*{4}}
\put(265,257){\circle*{4}}
\put(253,257){\circle*{4}}
\put(290,258){\circle*{4}}
\put(345,258){\circle*{4}}
\qbezier(290,258)(316,271)(345,258)
\qbezier(290,258)(319,247)(345,258)
\put(402,258){\circle*{4}}
\qbezier(345,258)(376,271)(402,258)
\qbezier(345,258)(375,246)(402,258)
\put(442,258){\circle*{4}}
\put(431,258){\circle*{4}}
\put(420,258){\circle*{4}}
\put(456,258){\circle*{4}}
\put(506,258){\circle*{4}}
\qbezier(456,258)(481,271)(506,258)
\qbezier(456,258)(482,247)(506,258)
\put(511,256){$v_{t}$}
\put(340,266){$v_{p}$}
\put(82,256){$\mathcal{D}$}
\put(258,220){$\mathcal{G}_{1}$}
\put(263,15){$\mathcal{G}_{2}$}
\put(121,137){$\mathbf{P}_{1}$}
\put(312,165){$\mathbf{P}_{2}$}
\put(363,111){$\mathbf{P}_{f}$}
\put(351,74){\circle*{4}}
\qbezier(289,75)(323,89)(351,74)
\qbezier(289,75)(323,61)(351,74)
\put(315,55){$\tilde{e}_{2}$}
\put(362,39){$\mathbf{P}_{0}$}
\put(112,256){$v_{0}$}
\put(150,242){$e_{1}$}
\put(181,265){$v_{1}$}
\put(234,264){$v_{2}$}
\put(204,242){$e_{2}$}
\put(477,242){$e_{t}$}
\put(274,278){$\mathbf{P}$}
\put(392,266){$v_{p+1}$}
\qbezier(400,74)(430,89)(463,74)
\qbezier(400,74)(432,61)(463,74)
\qbezier(463,74)(500,90)(526,75)
\qbezier(463,74)(502,62)(526,75)
\put(361,162){\circle*{4}}
\put(368,155){\circle*{4}}
\put(375,148){\circle*{4}}
\end{picture}
\end{center}

\begin{lemma} \label{le3,41}
$\mathcal{D}$ is a $k$-uniform hypergraph where $v_{0}\in V(\mathcal{D})$. Both $\mathbf{P}=v_{0}e_{1}v_{1}e_{2}v_{2}\cdots e_{t}v_{t}$ and
$\mathbf{P}_{0}=v_{0}\tilde{e}_{1}u_{1}\tilde{e}_{2}u_{2}\cdots \tilde{e}_{s}u_{s}$ are $k$-uniform hyperpaths where $\tilde{e}_{1}=\{v_{0}$, $v_{\varphi(1,1)}$, $v_{\varphi(1,2)}$, $\ldots$, $v_{\varphi(1,k-2)}$, $u_{1}\}$, $V(\mathcal{D})\cap V(\mathbf{P})=\{v_{0}\}$, $V(\mathcal{D})\cap V(\mathbf{P}_{0})=\{v_{0}\}$.
$\mathbf{P}_{1}$, $\mathbf{P}_{2}$, $\ldots$, $\mathbf{P}_{f}$ ($1\leq f\leq k-2$) are $k$-uniform hyperpaths attached respectively at vertices $v_{\varphi(1,1)}$, $v_{\varphi(1,2)}$, $\ldots$, $v_{\varphi(1,f)}$ in $\tilde{e}_{1}$ satisfying $1\leq L(\mathbf{P}_{1})\leq L(\mathbf{P}_{2})\leq\cdots\leq L(\mathbf{P}_{f})\leq  L(\mathbf{P}_{0})-1$, $\sum^{f}_{i=0}L(\mathbf{P}_{i})=t$, $V(\mathcal{D})\cap V(\mathbf{P}_{i})=\emptyset$ for $1\leq i\leq f$. $K$-uniform hypergraph $\mathcal{G}_{1}$ is a $\mathcal{G}(\mathcal{D}v_{0}; 0, t; v_{t})$ consisting of $\mathcal{D}$ and $\mathbf{P}$;
$k$-uniform hypergraph $\mathcal{G}_{2}$ consists of $\mathcal{D}$ and $\mathbf{P}_{0}$, $\mathbf{P}_{1}$, $\mathbf{P}_{2}$, $\ldots$, $\mathbf{P}_{f}$ (see Fig. 3.5).
Then $\rho(\mathcal{G}_{1})<\rho(\mathcal{G}_{2})$.
\end{lemma}

\begin{proof}
For brevity, we denote by $L_{i}=L(\mathbf{P}_{i})$ for $0\leq i\leq f$.  Without loss of generality, we suppose $f=3$ next.

In $\mathcal{G}_{1}$, denote by $e_{i}=\{v_{i-1}$, $v_{a(i,1)}$, $v_{a(i,2)}$, $\ldots$, $v_{a(i,k-2)}$, $v_{i}\}$ for $1\leq i\leq t-1$, $e_{t}=\{v_{t-1}$, $v_{a(t,1)}$, $v_{a(t,2)}$, $\ldots$, $v_{a(t,k-2)}$, $v_{a(t,k-1)}\}$ where $v_{a(t,k-1)}=v_{t}$. Assume that in $\mathcal{D}$, the edges incident with $v_{0}$ are $\varepsilon_{1}$, $\varepsilon_{2}$, $\ldots$, $\varepsilon_{\eta}$. Let $X$ be
the principal eigenvector of $\mathcal{G}_{1}$, and $x_{v_{p}}=\max\{v_{i}: 0\leq i\leq t\}$. By Lemma \ref{le3,40}, we know that $p\leq t-1$. By Lemma \ref{le3,08} and Lemma \ref{le3,08,02}, we know that $x_{v_{a(i,j)}}=x_{v_{a(i,z)}}<\min\{x_{v_{i-1}}, x_{v_{i}}\}$ for $2\leq j,z\leq k-2$ where $1\leq i\leq t-1$, $x_{v_{a(t,j)}}=x_{v_{a(t,z)}}$ for $2\leq j,z\leq k-1$. By Lemma \ref{le3,40}, we know that if $p>0$, then $x_{v_{i}}\leq x_{v_{i+1}}$ for $0\leq i\leq p-1$; if $p\geq0$, then $x_{v_{i}}\geq x_{v_{i+1}}$ for $p\leq i\leq t-1$.

{\bf Case 1} $p>0$.

{\bf Subcase 1.1} $t-p\geq L_{0}$ (see Fig. 3.6).
By Lemma \ref{le3,41,0}, there exists $1\leq j\leq3$ such that $t-p> L_{1}+\cdots +L_{j}$, but $L_{0}+L_{1}+\cdots +L_{j}> t-p$. Without loss of generality, we suppose $j=2$. Now $t-L_{1}-L_{2}> p$, $t-L_{0}-L_{1}-L_{2}< p$. Note that $L_{3}=t-L_{0}-L_{1}-L_{2}$. For brevity and convenience, we let $\xi=t-L_{1}$, $\eta=t-L_{1}-L_{2}$.

\setlength{\unitlength}{0.66pt}
\begin{center}
\begin{picture}(701,114)
\qbezier(12,45)(12,56)(21,64)\qbezier(21,64)(31,72)(45,72)\qbezier(45,72)(58,72)(68,64)\qbezier(68,64)(78,56)(78,45)\qbezier(78,45)(78,33)(68,25)
\qbezier(68,25)(58,18)(45,18)\qbezier(45,18)(31,18)(21,25)\qbezier(21,25)(12,33)(12,45)
\put(78,45){\circle*{4}}
\put(129,45){\circle*{4}}
\put(141,44){\circle*{4}}
\put(153,44){\circle*{4}}
\put(164,44){\circle*{4}}
\put(487,45){\circle*{4}}
\put(542,45){\circle*{4}}
\put(599,45){\circle*{4}}
\put(610,44){\circle*{4}}
\put(622,44){\circle*{4}}
\put(633,44){\circle*{4}}
\put(644,45){\circle*{4}}
\put(700,45){\circle*{4}}
\put(271,-9){Fig. 3.6. $\mathcal{G}_{1}\ (t-p\geq L_{0})$}
\put(32,41){$\mathcal{D}$}
\put(63,44){$v_{0}$}
\put(535,54){$v_{\xi}$}
\put(587,53){$v_{\xi+1}$}
\put(702,47){$v_{t}$}
\put(124,51){$v_{1}$}
\put(330,44){\circle*{4}}
\put(386,44){\circle*{4}}
\put(440,44){\circle*{4}}
\put(474,43){\circle*{4}}
\put(463,43){\circle*{4}}
\put(451,43){\circle*{4}}
\put(372,54){$v_{\eta-1}$}
\put(434,50){$v_{\eta}$}
\put(473,54){$v_{\xi-1}$}
\put(630,52){$v_{t-1}$}
\put(100,29){$e_{1}$}
\put(407,30){$e_{\eta}$}
\put(667,29){$e_{t}$}
\qbezier(78,45)(104,57)(129,45)
\qbezier(78,45)(103,34)(129,45)
\qbezier(330,44)(360,57)(386,44)
\qbezier(330,44)(360,33)(386,44)
\qbezier(386,44)(413,57)(440,44)
\qbezier(386,44)(413,35)(440,44)
\qbezier(487,45)(517,57)(542,45)
\qbezier(542,45)(570,57)(599,45)
\qbezier(644,45)(671,55)(700,45)
\qbezier(487,45)(518,33)(542,45)
\qbezier(542,45)(569,33)(599,45)
\qbezier(644,45)(672,34)(700,45)
\put(557,29){$e_{\xi+1}$}
\put(174,44){\circle*{4}}
\put(226,44){\circle*{4}}
\qbezier(174,44)(199,59)(226,44)
\qbezier(174,44)(198,33)(226,44)
\put(292,43){\circle*{4}}
\put(304,43){\circle*{4}}
\put(315,43){\circle*{4}}
\put(215,53){$v_{L_{3}}$}
\put(280,44){\circle*{4}}
\qbezier(226,44)(255,59)(280,44)
\qbezier(226,44)(255,31)(280,44)
\put(194,28){$e_{L_{3}}$}
\end{picture}
\end{center}

{\bf Subcase 1.1.1} $x_{v_{\eta}}\geq x_{v_{L_{3}}}$. By Lemma \ref{le3,40}, we know that $x_{v_{\eta-1}}\geq x_{v_{L_{3}+1}}$.

{\bf Subcase 1.1.1.1} $x_{v_{0}}\leq x_{v_{a(\eta,1)}}$, $x_{v_{L_{3}}}\leq x_{v_{a(\eta,2)}}$, $x_{v_{\xi}}\leq x_{v_{a(\eta,3)}}$. Let $\varepsilon^{'}_{i}=(\varepsilon_{i}\setminus \{v_{0}\})\cup \{v_{a(\eta,1)}\}$ for $1\leq i\leq\eta$, $e^{'}_{L_{3}}=(e_{L_{3}}\setminus\{v_{L_{3}}\})\cup\{v_{a(\eta,2)}\}$, $e^{'}_{\xi+1}=(e_{\xi+1}\setminus\{v_{\xi}\})\cup\{v_{a(\eta,3)}\}$;
$\mathcal {S}_{i}=\sum_{v\in (\varepsilon_{i}\setminus \{v_{0}\})}x_{v}$ for $i=1$, $2$, $\ldots$, $\eta$, $\mathbb{S}_{L_{3}}=\sum_{v\in e_{L_{3}}\setminus\{v_{L_{3}}\})}x_{v}$,
$\mathbb{S}_{\xi+1}=\sum_{v\in (e_{\xi+1}\setminus\{v_{\xi}\})}x_{v}$.

Let
$$\mathcal{G}^{'}_{1}=\mathcal{G}_{1}-\sum^{\eta}_{i=1}\varepsilon_{i}+\sum^{\eta}_{i=1}\varepsilon^{'}_{i}-e_{L_{3}}+e^{'}_{L_{3}}-e_{\xi+1}+e^{'}_{\xi+1}.$$ Then
$$X^{T}\mathcal {A}(\mathcal{G}^{'}_{1})X-X^{T}\mathcal {A}(\mathcal{G}_{1})X=\frac{2}{k-1}\{(x_{v_{a(\eta,1)}}-x_{v_{0}})\sum^{\eta}_{i=1}\mathcal {S}_{i}+(x_{v_{a(\eta,2)}}-x_{v_{L_{3}}})\mathbb{S}_{L_{3}}+(x_{v_{a(\eta,3)}}-x_{v_{\xi}})\mathbb{S}_{\xi+1}\}\geq 0.$$
It follows that $\rho(\mathcal{G}^{'}_{1})\geq \rho(\mathcal{G}_{1})$. Suppose $\rho(\mathcal{G}^{'}_{1})= \rho(\mathcal{G}_{1})$. Then $\rho(\mathcal{G}^{'}_{1})=X^{T}\mathcal {A}(\mathcal{G}^{'}_{1})X=X^{T}\mathcal {A}(\mathcal{G}_{1})X= \rho(\mathcal{G}_{1})$. Hence $X$ is also the principal eigenvector of $\mathcal{G}^{'}_{1}$ and $\mathcal {A}(\mathcal{G}^{'}_{1})X=\mathcal {A}(\mathcal{G}_{1})X$.
But a contradiction comes immediately because $(\mathcal {A}(\mathcal{G}^{'}_{1})X)_{v_{a(\eta,1)}}>(\mathcal {A}(\mathcal{G}_{1})X)_{v_{a(\eta,1)}}$.
Thus it follows that $\rho(\mathcal{G}^{'}_{1})> \rho(\mathcal{G}_{1})$. Note that $\mathcal{G}^{'}_{1}\cong \mathcal{G}_{2}$. Then $\rho(\mathcal{G}_{2})> \rho(\mathcal{G}_{1})$.

{\bf Subcase 1.1.1.2} $x_{v_{0}}> x_{v_{a(\eta,1)}}$, $x_{v_{L_{3}}}\leq x_{v_{a(\eta,2)}}$, $x_{v_{t-L_{1}}}\leq x_{v_{a(\eta,3)}}$.
Let $e^{'}_{1}=(e_{1}\setminus \{v_{0}\})\cup \{v_{a(\eta,1)}\}$, $e^{'}_{\eta}=(e_{\eta}\setminus\{v_{a(\eta,1)}\})\cup\{v_{0}\}$, $e^{'}_{L_{3}}=(e_{L_{3}}\setminus\{v_{L_{3}}\})\cup\{v_{a(\eta,2)}\}$, $e^{'}_{\xi+1}=(e_{\xi+1}\setminus\{v_{\xi}\})\cup\{v_{a(\eta,3)}\}$;
$\mathbb{S}_{1}=\sum_{v\in (e_{1}\setminus \{v_{0}\})}x_{v}$, $\mathbb{S}_{\eta}=\sum_{v\in (e_{\eta}\setminus \{v_{a(\eta,1)}\})}x_{v}$, $\mathbb{S}_{L_{3}}=\sum_{v\in e_{L_{3}}\setminus\{v_{L_{3}}\})}x_{v}$,
$\mathbb{S}_{\xi+1}=\sum_{v\in (e_{\xi+1}\setminus\{v_{\xi}\})}x_{v}$. Note $x_{v_{\eta}}\geq x_{v_{t-L_{0}-L_{1}-L_{2}}}\geq x_{v_{0}}$, $x_{\eta-1}\geq x_{t-L_{0}-L_{1}-L_{2}+1}=x_{L_{3}+1}$. Using Lemma \ref{le3,08}, we get $x_{v_{a(1,j)}}<x_{a(L_{3}+1,j)}<x_{v_{a(\eta,j)}}<\min\{x_{v_{\eta-1}}, x_{v_{\eta}}\}$ for $1\leq j\leq k-2$. As a result, it follows that $\mathbb{S}_{1}\leq \mathbb{S}_{\eta}$.

Let
$$\mathcal{G}^{'}_{1}=\mathcal{G}_{1}-e_{1}+e^{'}_{1}-e_{\eta}+e^{'}_{\eta}-e_{L_{3}}+e^{'}_{L_{3}}-e_{\xi+1}+e^{'}_{\xi+1}.$$ Then
$$X^{T}\mathcal {A}(\mathcal{G}^{'}_{1})X-X^{T}\mathcal {A}(G_{1})X=\frac{2}{k-1}\{(x_{v_{0}}-x_{v_{a(\eta,1)}})\mathbb{S}_{\eta}-(x_{v_{0}}-x_{v_{a(\eta,1)}})\mathbb{S}_{1}
+(x_{v_{a(\eta,2)}}-x_{v_{L_{3}}})\mathbb{S}_{L_{3}}+(x_{v_{a(\eta,3)}}-x_{v_{\xi}})\mathbb{S}_{\xi+1}\}$$$$\geq 0.\hspace{8.2cm}$$
Thus it follows that $\rho(\mathcal{G}^{'}_{1})\geq \rho(\mathcal{G}_{1})$. Note that $(\mathcal {A}(\mathcal{G}^{'}_{1})X)_{v_{\eta}}>(\mathcal {A}(\mathcal{G}_{1})X)_{v_{\eta}}$ and $\mathcal{G}^{'}_{1}\cong \mathcal{G}_{2}$. As Subcase 1.1.1.1, we get that $\rho(\mathcal{G}_{2})> \rho(\mathcal{G}_{1})$.

{\bf Subcase 1.1.1.3} $x_{v_{0}}\leq x_{v_{a(\eta,1)}}$, $x_{v_{L_{3}}}> x_{v_{a(\eta,2)}}$, $x_{v_{\xi}}\leq x_{v_{a(\eta, 3)}}$. Let $\varepsilon^{'}_{i}=(\varepsilon_{i}\setminus \{v_{0}\})\cup \{v_{a(\eta,1)}\}$ for $1\leq i\leq\eta$,
$e^{'}_{\eta}=(e_{\eta}\setminus\{v_{a(\eta,2)}\})\cup\{v_{L_{3}}\}$, $e^{'}_{L_{3}+1}=(e_{L_{3}+1}\setminus\{v_{L_{3}}\})\cup\{v_{a(\eta,2)}\}$, $e^{'}_{\xi+1}=(e_{\xi+1}\setminus\{v_{\xi}\})\cup\{v_{a(\eta,3)}\}$.

Let
$$G^{'}_{1}=G_{1}-\sum^{\eta}_{i=1}\varepsilon_{i}+\sum^{\eta}_{i=1}\varepsilon^{'}_{i}-e_{L_{3}+1}+e^{'}_{L_{3}+1}-e_{\eta}+e^{'}_{\eta}-e_{\xi+1}+e^{'}_{\xi+1}.$$
As Subcase 1.1.1.1 and Subcase 1.1.1.2, we get that $\rho(\mathcal{G}^{'}_{1})> \rho(\mathcal{G}_{1})$, and $\rho(\mathcal{G}_{2})> \rho(\mathcal{G}_{1})$.

In the same way, for the subcases: (i) $x_{v_{0}}\leq x_{v_{a(\eta,1)}}$, $x_{v_{L_{3}}}\leq x_{v_{a(\eta,2)}}$, $x_{v_{\xi}}> x_{v_{a(\eta,3)}}$; (ii) $x_{v_{0}}> x_{v_{a(\eta,1)}}$, $x_{v_{L_{3}}}> x_{v_{a(\eta,2)}}$, $x_{v_{\xi}}\leq x_{v_{a(\eta,3)}}$; (iii) $x_{v_{0}}> x_{v_{a(\eta,1)}}$, $x_{v_{L_{3}}}> x_{v_{a(\eta,2)}}$, $x_{v_{\xi}}> x_{v_{a(\eta,3)}}$; (iv) $x_{v_{0}}\leq x_{v_{a(\eta,1)}}$, $x_{v_{L_{3}}}> x_{v_{a(\eta,2)}}$, $x_{v_{\xi}}> x_{v_{a(\eta,3)}}$;  (v) $x_{v_{0}}> x_{v_{a(\eta,1)}}$, $x_{v_{L_{3}}}\leq x_{v_{a(\eta,2)}}$, $x_{v_{\xi}}> x_{v_{a(\eta,3)}}$, we get that $\rho(\mathcal{G}_{2})> \rho(\mathcal{G}_{1})$.

As a result, from the above narrations for Subcase 1.1.1 that $x_{v_{\eta}}\geq x_{v_{L_{3}}}$, we get $\rho(\mathcal{G}_{2})> \rho(\mathcal{G}_{1})$.

{\bf Subcase 1.1.2} $x_{v_{\eta}}< x_{v_{L_{3}}}$. By Lemma \ref{le3,40}, we know that $x_{v_{\eta-1}}\leq x_{v_{L_{3}+1}}$. By considering the comparisons between $x_{v_{0}}$ and $x_{v_{a(L_{3}+1,1)}}$, between $x_{v_{\xi}}$ and $x_{v_{a(L_{3}+1,2)}}$, between $x_{v_{\eta}}$ and $x_{v_{a(L_{3}+1,3)}}$, as Subcase 1.1.1, we get that $\rho(\mathcal{G}_{2})> \rho(\mathcal{G}_{1})$.

{\bf Subcase 1.2} $t-p< L_{0}$ (see Fig. 3.7). Let $\omega= L_{1}+L_{2}$, $\varphi=L_{1}+L_{2}+L_{3}+1$. By considering the comparisons between $x_{v_{0}}$ and $x_{v_{a(\varphi,1)}}$, between $x_{v_{L_{1}}}$ and $x_{v_{a(\varphi,2)}}$, between $x_{v_{\omega}}$ and $x_{v_{a(\varphi,3)}}$, as Subcase 1.1, we get that $\rho(\mathcal{G}_{2})> \rho(\mathcal{G}_{1})$.

\setlength{\unitlength}{0.66pt}
\begin{center}
\begin{picture}(713,183)
\qbezier(0,45)(0,56)(9,64)\qbezier(9,64)(19,72)(33,72)\qbezier(33,72)(46,72)(56,64)\qbezier(56,64)(66,56)(66,45)\qbezier(66,45)(66,33)(56,25)
\qbezier(56,25)(46,18)(33,18)\qbezier(33,18)(19,18)(9,25)\qbezier(9,25)(0,33)(0,45)
\put(66,45){\circle*{4}}
\put(117,45){\circle*{4}}
\put(129,44){\circle*{4}}
\put(141,44){\circle*{4}}
\put(152,44){\circle*{4}}
\put(475,45){\circle*{4}}
\put(546,45){\circle*{4}}
\put(611,45){\circle*{4}}
\put(622,44){\circle*{4}}
\put(634,44){\circle*{4}}
\put(645,44){\circle*{4}}
\put(656,45){\circle*{4}}
\put(712,45){\circle*{4}}
\put(292,-9){Fig. 3.8. $\mathcal{G}_{1}\ (p= 0)$}
\put(20,44){$D$}
\put(49,43){$v_{0}$}
\put(537,53){$v_{\varsigma}$}
\put(599,53){$v_{\varsigma+1}$}
\put(716,44){$v_{t}$}
\put(110,52){$v_{1}$}
\put(318,44){\circle*{4}}
\put(374,44){\circle*{4}}
\put(428,44){\circle*{4}}
\put(462,43){\circle*{4}}
\put(451,43){\circle*{4}}
\put(439,43){\circle*{4}}
\put(366,52){$v_{\kappa}$}
\put(413,52){$v_{\kappa+1}$}
\put(464,53){$v_{\varsigma-1}$}
\put(88,29){$e_{1}$}
\put(387,30){$e_{\kappa+1}$}
\put(680,30){$e_{t}$}
\qbezier(66,45)(92,57)(117,45)
\qbezier(66,45)(91,34)(117,45)
\qbezier(318,44)(348,57)(374,44)
\qbezier(318,44)(348,33)(374,44)
\qbezier(374,44)(401,57)(428,44)
\qbezier(374,44)(401,35)(428,44)
\qbezier(475,45)(521,57)(546,45)
\qbezier(546,45)(574,57)(611,45)
\qbezier(656,45)(683,55)(712,45)
\qbezier(475,45)(522,33)(546,45)
\qbezier(546,45)(573,33)(611,45)
\qbezier(656,45)(684,34)(712,45)
\put(567,28){$e_{\varsigma+1}$}
\put(162,44){\circle*{4}}
\put(214,44){\circle*{4}}
\qbezier(162,44)(187,59)(214,44)
\qbezier(162,44)(186,33)(214,44)
\put(280,43){\circle*{4}}
\put(292,43){\circle*{4}}
\put(303,43){\circle*{4}}
\put(204,53){$v_{L_{0}}$}
\put(268,44){\circle*{4}}
\qbezier(214,44)(243,59)(268,44)
\qbezier(214,44)(243,31)(268,44)
\put(182,26){$e_{L_{0}}$}
\qbezier(24,81)(25,81)(24,81)
\qbezier(24,81)(24,81)(25,81)
\qbezier(0,159)(0,168)(9,175)\qbezier(9,175)(18,182)(31,182)\qbezier(31,182)(43,182)(52,175)\qbezier(52,175)(62,168)(62,159)\qbezier(62,159)(62,149)(52,142)
\qbezier(52,142)(43,135)(31,135)\qbezier(31,136)(18,136)(9,142)\qbezier(9,142)(0,149)(0,159)
\put(61,159){\circle*{4}}
\put(114,159){\circle*{4}}
\qbezier(61,159)(89,175)(114,159)
\qbezier(61,159)(86,147)(114,159)
\put(127,158){\circle*{4}}
\put(139,158){\circle*{4}}
\put(150,158){\circle*{4}}
\put(163,158){\circle*{4}}
\put(214,158){\circle*{4}}
\qbezier(163,158)(188,174)(214,158)
\qbezier(163,158)(188,147)(214,158)
\put(265,158){\circle*{4}}
\qbezier(214,158)(240,173)(265,158)
\qbezier(214,158)(241,149)(265,158)
\put(278,157){\circle*{4}}
\put(290,157){\circle*{4}}
\put(301,157){\circle*{4}}
\put(314,158){\circle*{4}}
\put(368,158){\circle*{4}}
\qbezier(314,158)(345,173)(368,158)
\qbezier(314,158)(345,147)(368,158)
\put(425,159){\circle*{4}}
\qbezier(368,158)(398,174)(425,159)
\qbezier(368,158)(399,147)(425,159)
\put(439,159){\circle*{4}}
\put(451,159){\circle*{4}}
\put(462,159){\circle*{4}}
\put(475,159){\circle*{4}}
\put(541,159){\circle*{4}}
\qbezier(475,159)(507,175)(541,159)
\qbezier(475,159)(507,148)(541,159)
\put(602,159){\circle*{4}}
\qbezier(541,159)(577,174)(602,159)
\qbezier(541,159)(575,148)(602,159)
\put(617,159){\circle*{4}}
\put(629,159){\circle*{4}}
\put(640,159){\circle*{4}}
\put(653,160){\circle*{4}}
\put(708,160){\circle*{4}}
\qbezier(653,160)(678,176)(708,160)
\qbezier(653,160)(678,148)(708,160)
\put(19,158){$D$}
\put(44,158){$v_{0}$}
\put(110,164){$v_{1}$}
\put(81,142){$e_{1}$}
\put(204,168){$v_{L_{1}}$}
\put(179,143){$e_{L_{1}}$}
\put(360,167){$v_{\omega}$}
\put(336,142){$e_{\omega}$}
\put(532,168){$v_{\varphi}$}
\put(593,168){$v_{\varphi+1}$}
\put(559,144){$e_{\varphi+1}$}
\put(712,158){$v_{t}$}
\put(676,143){$e_{t}$}
\put(293,102){Fig. 3.7. $\mathcal{G}_{1}\ (t-p< L_{0})$}
\put(412,168){$v_{\omega+1}$}
\end{picture}
\end{center}

{\bf Case 2} $p=0$ (see Fig. 3.8). Let $\kappa=L_{0}+L_{1}$, $\varsigma=L_{0}+L_{1}+L_{2}$. By considering the comparisons between $x_{v_{a(1,1)}}$ and $x_{v_{L_{0}}}$, between $x_{v_{a(1,2)}}$ and $x_{v_{\kappa}}$,
between $x_{v_{a(1,3)}}$ and $x_{v_{\varsigma}}$, as Case 1, we get that $\rho(\mathcal{G}_{2})> \rho(\mathcal{G}_{1})$.

This completes the proof. \ \ \ \ \ $\Box$
\end{proof}

\begin{Prof}
For $k$-uniform supertrees of order $n$, using Lemma \ref{le3,39} and Lemma \ref{le3,41} repeatedly gets the result.
This completes the proof. \ \ \ \ \ $\Box$
\end{Prof}

\begin{Proff}
Note that $\mathcal{S}^{\ast}(n, k)$ is the $k$th power of the ordinary star $S^{\ast}(\frac{n-1}{k-1}+1, 2)$, $\mathcal{P}(n, k)$ is the $k$th power of the ordinary path $P(\frac{n-1}{k-1}+1, 2)$. Then the result follows from Theorem \ref{th01.01} and Theorem \ref{th01.03}.
This completes the proof. \ \ \ \ \ $\Box$
\end{Proff}

\small {

}

\end{document}